\documentclass[12pt, reqno, twoside]{amsart}

\usepackage[all]{xy}

\usepackage{euscript}

\usepackage[totalwidth=489.05pt, totalheight=722.70pt]{geometry}

\usepackage{times}

\usepackage{hyperref}

\hypersetup{pdfauthor={N. Raghavendra}}%
\hypersetup{pdftitle={Binary Trees and Fibred Categories}}%
\hypersetup{pdfsubject={Abstract formalism for binary trees}}%
\hypersetup{pdfkeywords={Binary tree, fibred category, binary graph,
    binary automaton, transition system}}
\hypersetup{pdfpagemode=true}
\hypersetup{colorlinks=false}%
\hypersetup{urlcolor=blue}%
\hypersetup{citecolor=red}%

\theoremstyle{definition}%
\newtheorem{definition}{Definition}[section]%
\newtheorem{notation}[definition]{Notation}%

\theoremstyle{remark}%
\newtheorem{remark}[definition]{Remark}%
\newtheorem{example}[definition]{Example}%

\theoremstyle{plain}%
\newtheorem{theorem}[definition]{Theorem}%
\newtheorem{proposition}[definition]{Proposition}%
\newtheorem{corollary}[definition]{Corollary}%

\renewcommand{\mathcal}{\EuScript}

\newcommand{\N}{\ensuremath{\mathbf{N}}}

\newcommand{\emptyword}{\ensuremath{\epsilon}}
\newcommand{\set}[1]{\ensuremath{\{ #1 \}}}

\newcommand{\length}[1]{\ensuremath{\vert #1 \vert}}
\newcommand{\suchthat}{\ensuremath{\, \vert \,}}
\newcommand{\Ob}[1]{\ensuremath{\operatorname{Ob} (#1)}}
\newcommand{\Hom}[3][]{\ensuremath{\operatorname{Hom}_{#1} (#2, #3)}}
\DeclareMathOperator{\HomCatOperator}{\ensuremath{\mathbf{Hom}}}
\newcommand{\HomCat}[3][]{%
  \ensuremath{\HomCatOperator_{#1} (#2, #3)}}
\newcommand{\Cat}{\ensuremath{\mathbf{Cat}}}
\newcommand{\id}[1]{\ensuremath{\mathbf{1}_{#1}}}
\DeclareMathOperator{\SplitOperator}{\ensuremath{\mathbf{Split}}}
\newcommand{\Split}[1]{\ensuremath{\SplitOperator_{#1}}}
\newcommand{\restrict}[2]{\ensuremath{#1 \vert_{#2}}}
\newcommand{\powerset}[1]{\ensuremath{\mathbf{2}^{#1}}}
\newcommand{\powersetcat}[1]{\ensuremath{\mathcal{P}_{#1}}}

\begin{document}

\title[Formalism for binary trees]{Binary trees and fibred categories}

\author{N.~Raghavendra}


\address{Advanced Technology Centre, Tata Consultancy Services, 5-9-62
  (6th Floor), Khan Lateef Khan Estate, Fateh Maidan Road, Hyderabad
  500 001, India}%

\email{raghu@atc.tcs.co.in}



\begin{abstract}
  We develop a purely set-theoretic formalism for binary trees and
  binary graphs. We define a category of binary automata, and display
  it as a fibred category over the category of binary graphs. We also
  relate the notion of binary graphs to transition systems, which
  arise in the theory of concurrent computing.
\end{abstract}


\maketitle

\tableofcontents

\setlength{\parskip}{0.45\baselineskip}

\section{Introduction}
\label{sec:introduction}

This work arose out of certain issues which we encountered while
trying to understand certain well-known algorithms in computational
biology, whose aim is to build trees out of various data. We will
discuss these issues presently. In any case, they motivated us to
develop an abstract formalism for binary trees, which is the subject
of this article. We found that our formalism fits neatly into the
framework of fibred categories, which were originally studied in
algebraic geometry. Moreover, it has connections with the notion of
transition systems which arise in the theory of concurrent
computation.

The problems which motivated this work are as follows. First, there
are several statements in tree-construction algorithms which assert
the uniqueness of some tree. One can ask certain natural questions in
such a situation. To begin with, what does one mean by uniqueness: for
instance, does one mean uniqueness in the sense of equality of trees,
or does one mean uniqueness up to an isomorphism?  Next, what is the
precise property of trees, which makes the tree in question unique?
Can one formulate this property as a universal property in a suitable
category? Thus, we were led to seek an abstract framework for binary
trees.

Secondly, the standard descriptions of trees are either recursive, as
in \cite[Section 2.3]{Knuth:1997:FA}, or are based on graph theory, as
in \cite[Section 2.4]{Aho:1974:DAC}. We feel that these descriptions
have certain drawbacks. On the one hand, a recursive definition of
trees sets them apart from other important classes of mathematical
objects such as groups, finite automata etc., which are described in
terms of sets and maps between them, whereas a set-theoretic
description of trees would bring them into a general framework,
consistent with the definitions of several other mathematical
structures. On the other hand, a treatment of trees based on graph
theory entails carrying, to some extent, the formalism of the latter
subject, which can be avoided if we, \textit{ab initio}, develop a
definition of trees which does not rely on graph theory.

In view of the above discussion, in this article we develop a purely
set-theoretic treatment of binary trees. We will also discuss the
relationship between our definitions, and the standard graph-theoretic
notions regarding trees.

We now sketch the contents of the various sections. In Section
\ref{sec:binarygraphs}, we develop a set-theoretic formalism for
binary graphs. In Section \ref{sec:bitstrings}, we establish our
notation for the free monoid on the set of bits, which monoid acts on
the set of nodes of any binary graph, as described in Section
\ref{sec:action}. Transition diagrams, defined in Section
\ref{sec:transitiondiagrams}, establish a connection between our model
for binary graphs and their standard graph-theoretic definition. In
Section \ref{sec:binaryforests}, we describe binary trees as a
specialization of binary forests. We introduce the category of binary
automata in Section \ref{sec:binaryautomata}, and we relate binary
graphs to transition systems. The important notion of a fibred
category is recalled in Section \ref{sec:fibredcategories}, and in the
two subsequent Sections, we recall a criterion for a category to be
fibred, and the notion of split categories. To assist the reader, and
to make this paper self-contained, we give complete definitions and
proofs regarding fibred categories, even though the material may be
well-known to specialists in certain areas of algebraic geometry or
category theory. In Section \ref{sec:binaryfibredcat}, we display the
category of binary automata as a fibred category over the category of
binary graphs. In Section \ref{sec:conclusion}, we summarise the
article, and mention some ideas for future work.

\section{Binary graphs}
\label{sec:binarygraphs}

We begin our discussion of binary trees with the definition of a more
general object, which we call a \emph{binary graph}. This helps us to
introduce several constructions which are useful in our later
discussion of binary trees.

\begin{definition}
  \label{def:binary graph}
  A \emph{binary graph} is a pair
  $\Gamma = (Q, \delta)$,
  where
  \begin{enumerate}
  \item $Q$ is a finite set, whose elements are called \emph{nodes}.
  \item $\delta: Q \times \Sigma \rightarrow Q$ is a function, called
    the \emph{transition function}, where $\Sigma = \set{0,1}$.
  \end{enumerate}
\end{definition}

Given such a binary graph $\Gamma$, we get a pair of functions
$\lambda : Q \rightarrow Q$ and $\rho : Q \rightarrow Q$, defined by
\begin{displaymath}
  \lambda (q) = \delta (q, 0), \quad \text{and} \quad
  \rho (q) = \delta (q, 1), \quad \text{for all $q \in Q$}.
\end{displaymath}
We call $\lambda (q)$ the \emph{left child} of $q$, and $\rho (q)$ the
\emph{right child} of $q$.

\begin{definition}
  \label{def:morphismbinarygraph}
  Let $\Gamma = (Q, \delta)$ and $\Gamma' = (Q', \delta')$ be binary
  graphs. A \emph{morphism} from $\Gamma$ to $\Gamma'$ is a function
  $f : Q \rightarrow Q'$, such that
  $f (\delta (q, a)) =\delta' (f (q), a)$, for all $q \in Q$ and for
  all letters $a \in \Sigma$. We say that $\Gamma$ is a
  \emph{binary subgraph} of $\Gamma'$ if $Q$ is a subset of $Q'$, and
  if the inclusion $i : Q \hookrightarrow Q'$ is a morphism of binary
  graphs.
\end{definition}

\begin{example}
  \label{exa:binarygraph1}
  Let $Q = \set{q_0, q_1}$, and define
  $\delta : Q \times \Sigma \rightarrow Q$ by $\delta (q, 0) = q_0$
  and $\delta (q, 1) = q_1$, for all $q \in Q$. It is convenient to
  define the transition function $\delta$ through a table, whose rows
  are indexed by the letters in $\Sigma$, whose columns are indexed by
  the nodes in $Q$, and whose $(a, q)$-th entry equals
  $\delta (q, a)$. In this notation, the above function is defined by
  the following table:
  \begin{center}
    \begin{tabular}{|c|c|c|}
      \hline
      & $q_0$ & $q_1$ \\
      \hline
      0 & $q_0$ & $q_0$ \\
      1 & $q_1$ & $q_1$ \\
      \hline
    \end{tabular}
  \end{center}
\end{example}

\begin{example}
  \label{exa:binarygraph2}
  Let $Q = \set{q_0, q_1, q_2, q_3, q_4, q_5, q_6}$, and define
  $\delta : Q \times \Sigma \rightarrow Q$ by
  the following table:
  \begin{center}
    \begin{tabular}{|c|c|c|c|c|c|c|c|}
      \hline
      & $q_0$ & $q_1$ & $q_2$ & $q_3$ & $q_4$ & $q_5$ & $q_6$  \\
      \hline
      0 & $q_1$ & $q_3$ & $q_4$ & $q_3$ & $q_4$ & $q_6$ & $q_6$ \\
      1 & $q_2$ & $q_4$ & $q_5$ & $q_3$ & $q_4$ & $q_5$ & $q_6$ \\
      \hline
    \end{tabular}
  \end{center}
\end{example}

\begin{example}
  \label{exa:binarygraph3}
  Let $Q = \set{q_0, q_1, q_2, q_3, q_4, q_5, q_6, q_7}$, and define
  $\delta : Q \times \Sigma \rightarrow Q$ by
  the following table:
  \begin{center}
    \begin{tabular}{|c|c|c|c|c|c|c|c|c|}
      \hline
      & $q_0$ & $q_1$ & $q_2$ & $q_3$ & $q_4$ & $q_5$ & $q_6$ & $q_7$
      \\
      \hline
      0 & $q_1$ & $q_3$ & $q_2$ & $q_3$ & $q_4$ & $q_6$ & $q_6$ &
      $q_7$ \\
      1 & $q_2$ & $q_4$ & $q_2$ & $q_3$ & $q_4$ & $q_7$ & $q_6$ &
      $q_7$ \\
      \hline
    \end{tabular}
  \end{center}
\end{example}

\begin{example}
  \label{exa:binarygraph4}
  Let $Q = \set{q_0, q_1, q_2, q_3}$, and define
  $\delta : Q \times \Sigma \rightarrow Q$ by
  the following table:
  \begin{center}
    \begin{tabular}{|c|c|c|c|c|}
      \hline
      & $q_0$ & $q_1$ & $q_2$ & $q_3$ \\
      \hline
      0 & $q_1$ & $q_3$ & $q_2$ & $q_3$ \\
      1 & $q_2$ & $q_1$ & $q_2$ & $q_3$ \\
      \hline
    \end{tabular}
  \end{center}
\end{example}

\section{Bit strings}
\label{sec:bitstrings}

In this section, we establish our notation for the languages on the
set of bits. Let $\Sigma^{*}$ be the free monoid on the set $\Sigma$
which, we recall, is the set $\set{0,1}$. Thus, $\Sigma^{*}$ is a
disjoint union,
\begin{displaymath}
  \Sigma^{*} = \coprod_{n=0}^{\infty} \Sigma^{n},
\end{displaymath}
where $\Sigma^0$ is a singleton consisting of a symbol $\emptyword$,
and $\Sigma^n = \Sigma \times \cdots \times \Sigma$ ($n$ times). We
call elements of $\Sigma^{*}$ \emph{words} or \emph{bit strings}. In
particular, the word $\emptyword$ is called the \emph{empty word}. We
call the set $\Sigma$ an \emph{alphabet}, and refer to the elements of
$\Sigma$ as \emph{letters}. If a word
$w \in \Sigma^n$, we say that $w$ has length $n$, and write
$\length{w} = n$. Following standard practice, we denote a word
$w = (a_1, \ldots, a_n) \in \Sigma^n$ by $w = a_1 \cdots a_n$.There is
a natural multiplication law in the set $\Sigma^{*}$,
\begin{align*}
  \Sigma^{*} \times \Sigma^{*} & \rightarrow \Sigma^{*}, \\
  (v, w) & \mapsto vw,
\end{align*}
which is given by concatenation of words. Explicitly, we define
\begin{displaymath}
  \emptyword w = w \emptyword = w \quad
  \text{for all $w \in \Sigma^{*}$},
\end{displaymath}
and if $v = a_1 \cdots a_m$ and $w = b_1 \cdots b_n$, where
$a_i, b_j \in \Sigma$ for all $i = 1, \ldots, m$ and
$j = 1, \ldots, n$, then we set
\begin{displaymath}
  vw = a_1 \cdots a_m b_1 \cdots b_n.
\end{displaymath}
With this multiplication law, the set $\Sigma^{*}$ becomes a monoid,
whose identity element is the empty word $\emptyword$. It is clear
that $\length{vw} = \length{v} + \length{w}$ for all
$v, w \in \Sigma^{*}$. The monoid $\Sigma^{*}$ is free on the set
$\Sigma$, i.e., it has the following universal property: given any
monoid $M$, and given any function
$f : \Sigma \rightarrow M$, there exists a unique homomorphism of
monoids, $\tilde{f}: \Sigma^{*} \rightarrow M$, which extends $f$,
i.e., which makes the diagram
\begin{displaymath}
  \xymatrix{
    \Sigma \ar[dr]_{f} \ar[r]^i & \Sigma^{*} \ar[d]^{\tilde{f}} \\
    & M 
  }
\end{displaymath}
commute, where $i: \Sigma \rightarrow \Sigma^{*}$ is the natural
inclusion map.

\section{Action of bit strings}
\label{sec:action}

Let $\Gamma = (Q, \delta)$ be a binary graph. The universal property
of $\Sigma^{*}$ implies that there is a unique right action of the
monoid $\Sigma^{*}$ on $Q$,
\begin{align*}
  \hat{\delta} : Q \times \Sigma^{*} & \rightarrow Q, \\
  (q, w) & \mapsto qw,
\end{align*}
such that
\begin{displaymath}
  \hat{\delta} (q, aw) = \hat{\delta} (\delta (q, a), w) \quad
  \text{for all
    $q \in Q$, $a \in \Sigma$, and $w \in \Sigma^{*}$}.
\end{displaymath}
By abuse of notation, we denote $\hat{\delta}$ also by $\delta$.

For any node $q \in Q$, let
\begin{displaymath}
  \Sigma^{*}_{q} = \set{w \in \Sigma^{*} \suchthat qw = q}
\end{displaymath}
be the isotropy submonoid of $\Sigma^{*}$ at $q$. We say that a node
$q$ is a \emph{leaf} if it is fixed by $\Sigma^{*}$, i.e., if
$\Sigma^{*}_{q} = \Sigma^{*}$. If $q$ is not a leaf, i.e., if
$\Sigma^{*}_{q}$ is a proper submonoid of $\Sigma^{*}$, then we call
$q$ an \emph{internal node} of the binary graph $\Gamma$. Note that
$q$ is a leaf if and only if $\delta (q, a) = q$ for all
letters $a \in \Sigma$.

\begin{example}
  \label{exa:leaves}
  The binary graph of Example \ref{exa:binarygraph1}\ has no leaves.
  In Example \ref{exa:binarygraph2}, the nodes $q_3$, $q_4$ and $q_6$
  are the leaves. The leaves of $\Gamma$, in Example
  \ref{exa:binarygraph3}, are the nodes $q_2$, $q_3$, $q_4$, $q_6$ and
  $q_7$.
\end{example}

We say that an element $q$ of $Q$ is a \emph{parent} of $q'$ if
$q \neq q'$, and $\delta (q, a) = q'$ for some $a \in \Sigma$. In this
case $q$ has to be an internal node.  A \emph{path} in $\Gamma$ is a
pair $P = (q_0, w)$, where $q_0 \in Q$ is a node, and
$w \in \Sigma^{*}$ is a word, satisfying one of the following two
conditions:
\begin{enumerate}
\item Either $w = \emptyword$ is the empty word. Or,
\item $w = a_1 \cdots a_n$, with $a_i \in \Sigma$, and
  $q_i \neq q_{i-1}$, where $q_i = q_0 a_1 \cdots a_i$, for each
  $i = 1, \ldots, n$. (In particular, $q_{i-1}$ is a parent of $q_i$,
  for each $i = 1, \ldots, n$.)
\end{enumerate}
We say that this path \emph{starts} at $q_0$, \emph{ends} at
$q_0 w$, and has \emph{length} $\length{w}$.  A \emph{trivial path} is
one which has length $0$; it has to be of the form $(q_0,
\emptyword)$, for some node $q_0$. We say that $q$ is an
\emph{ancestor} of $q'$, or that $q'$ is a \emph{descendant} of $q$,
if there is a path starting at $q$ and ending at $q'$. If, in
addition, $q \neq q'$, then we say that $q$ is a
\emph{proper ancestor} of $q'$, or that $q'$ is a
\emph{proper descendant} of $q$.

\begin{definition}
  \label{def:InheritancePreorder}
  Let $\Gamma = (Q, \delta)$ be a binary graph. For any pair of nodes,
  $q, q' \in Q$, we write $q' \leq q$ (respectively, $q' < q$) if $q$
  is an ancestor (respectively, a proper ancestor) of $q'$. We write
  $q \geq q'$ (respectively, $q > q'$) if $q' \leq q$ (respectively,
  $q' < q$).  The relation $\leq$ on $Q$ is a preorder, i.e., it has
  all the properties of a partial order, except that of anti-symmetry.
  In other words, the relation $\leq$ is reflexive and transitive, but
  there may exist elements $q, q' \in Q$ such that
  $q \leq q'$ and $q' \leq q$, and yet $q \neq q'$. We call this
  preorder the \emph{inheritance preorder} on $Q$.
\end{definition}

\begin{remark}
  \label{rem:LeafPath}
  Any path which starts at a leaf is a trivial path.
\end{remark}

\begin{remark}
  \label{rem:morphismbinarygraph}
  Let $\Gamma = (Q, \delta)$ and $\Gamma' = (Q', \delta')$ be binary
  graphs. A morphism from $\Gamma$ to $\Gamma'$, introduced in
  Definition \ref{def:morphismbinarygraph}, is just a function
  $f : Q \rightarrow Q'$, which is $\Sigma^{*}$-equivariant, i.e.,
  which satisfies the condition that $f(qw) = f(q)w$ for all $q \in Q$
  and $w \in \Sigma^{*}$.  Similarly, if $\Gamma' = (Q', \delta')$ is
  a binary graph, and $Q$ is a $\Sigma^{*}$-invariant subset of $Q'$,
  which means that $qw \in Q$ for all $q \in Q$ and for all
  $w \in \Sigma^{*}$, we obtain a binary subgraph
  $\Gamma = (Q, \delta)$ of $\Gamma'$, where
  $\delta : Q \times \Sigma \rightarrow Q$ is the restriction of
  $\delta' : Q' \times \Sigma \rightarrow Q'$. Thus, binary subgraphs
  of $\Gamma'$ are in bijective correspondence with
  $\Sigma^{*}$-invariant subsets of $Q'$. Note that a subset $Q$ of
  $Q'$ is $\Sigma^{*}$-invariant if and only if
  $\delta' (q, a) \in Q$, for all $q \in Q$ and for all letters
  $a \in \Sigma$.
\end{remark}

For every node $q \in Q$, let
\begin{displaymath}
  \Sigma^{*} (q) = \set{qw \suchthat w \in \Sigma^{*}}
\end{displaymath}
be the orbit of $q$.

\begin{proposition}
  \label{pro:ancestor}
  In any binary graph $\Gamma = (Q, \delta)$, a node $q$ is an
  ancestor of $q'$ if and only if $q' \in \Sigma^{*} (q)$.
\end{proposition}

Before proving the Proposition, let us introduce a device which will
be useful in the proof. Let $\Gamma = (Q, \delta)$ be a binary graph,
and let $w \in \Sigma^{*}$. For each node $q_0 \in Q$, we shall define
a new word $w [q_0]$ as follows. If $w = \emptyword$, then set
$w [q_0] = w = \emptyword$. If $w$ is nonempty, write
$w = a_1 \cdots a_n$, where $a_i \in \Sigma$, and let
$q_i = q_0 a_1 \cdots a_i$, for $i = 1, \ldots, n$. If $q_i = q_{i-1}$
for all $i = 1, \ldots, n$, define $w [q_0] = \emptyword$. If not all
$q_i$ are equal, let $i_1 < \cdots < i_k$ be the elements of the set
\begin{displaymath}
  \set{i \suchthat 1 \leq i \leq n \; \text{and} \; q_i \neq q_{i-1}}.
\end{displaymath}
Let $b_r = a_{i_r}$, for $r = 1, \ldots, k$, and define
$w [q_0] = b_1 \cdots b_k$. The words $w$, $w [q_0]$ and the indices
$i_r$ have the following properties:
\begin{enumerate}
\item \label{item:ContractPath} The pair $P [q_0, w] = (q_0, w [q_0])$
  is a path, starting at $q_0$ and ending at $q_0 w$.
\item \label{item:ContractNodes} We have
  \begin{align*}
    q_i & = q_0, \quad \text{if $0 \leq i < i_1$}, \\
    q_i & = q_{i_r}, \quad \text{%
      if $1 \leq r < k$ and $i_r \leq i < i_{r+1}$, and} \\
    q_i & = q_{i_k}, \quad \text{if $i_k \leq i \leq n$.}
  \end{align*}
\item \label{item:ContractAi} The letters $a_i$ satisfy the following
  conditions:
  \begin{align*}
    q_0 a_i & = q_0, \quad \text{if $1 \leq i < i_1$}, \\
    q_{i_r} a_i & = q_{i_r}, \quad \text{%
      if $1 \leq r < k$ and $i_r < i < i_{r+1}$, and} \\
    q_{i_k} a_i & = q_{i_k}, \quad \text{if $i_k < i \leq n$.}
  \end{align*}
\item \label{item:ContractBi} The $b_i$ satisfy the following
  equations:
  \begin{align*}
    q_0 b_1 & = q_{i_1}, \\
    q_{i_r} b_{r+1} & = q_{i_{r+1}}, \quad \text{%
      if $1 \leq r < k$.}
  \end{align*}
\end{enumerate}

\begin{definition}
  \label{def:contraction}
  Let $\Gamma = (Q, \delta)$ be a binary graph, and let
  $w \in \Sigma^{*}$ be a word. Then, for each node $q_0 \in Q$, we
  call the word $w [q_0]$ defined above, the \emph{contraction} of $w$
  by $q_0$.
\end{definition}

\noindent \emph{Proof of Proposition \ref{pro:ancestor}}. Clearly, if
there is a path from $q$ to $q'$, then $q' \in \Sigma^{*} (q)$.
Conversely, suppose that $q' \in \Sigma^{*} (q)$, say $q' = qw$. Then,
by the device introduced above, there exists a word $w [q]$ such that
the pair $P = (q, w [q])$ is a path from $q$ to $q w = q'$. \qed

\begin{remark}
  \label{rem:OrbitGraph}
  If $q$ is a leaf in a binary graph $\Gamma = (Q, \delta)$, then
  $\Sigma^{*}$ fixes $q$, so we get a binary subgraph
  $\Gamma_q = (\set{q}, \delta_q)$, where $\delta_q$ is the unique map
  $\set{q} \times \Sigma \rightarrow \set{q}$. More generally, if $q$
  is any node, then we obtain a binary subgraph
  $\Gamma_q = (\Sigma^{*} (q), \delta_q)$, where $\Sigma^{*} (q)$ is
  the orbit of $q$, and
  $\delta_q : \Sigma^{*} (q) \times \Sigma \rightarrow \Sigma^{*} (q)$
  is the restriction of $\delta : Q \times \Sigma \rightarrow Q$. We
  call $\Gamma_q$ the \emph{orbit graph} of $q$. It is obvious that
  $\Gamma_q$ is a minimal binary subgraph of $\Gamma$, whenever $q$ is
  a leaf of $\Gamma$.
\end{remark}

The next proposition provides criteria , in terms of the action of
$\Sigma^{*}$, for the anti-symmetry of the inheritance preorder
$\leq$, defined above (see Definition \ref{def:InheritancePreorder}).
We say that a word $v \in \Sigma^{*}$ is a \emph{prefix} of another
word $w$, if $w = vz$ for some word $z \in \Sigma^{*}$. In this case
we say that the word $z$ is a \emph{suffix} of $w$. We say that $v$ is
a \emph{proper prefix} (respectively, \emph{proper suffix}) of $w$ if
$v$ is a prefix (respectively, suffix) of $w$,
$v \neq \emptyword$, and $v \neq w$. A subset $L$ of $\Sigma^{*}$ is
said to be \emph{prefix-closed}, if whenever $w \in L$, every prefix
$v$ of $w$ also belongs to $L$. The notion of a \emph{suffix-closed}
subset of $\Sigma^{*}$ is defined in a similar manner. Clearly, every
nonempty prefix-closed set has to contain the empty word $\emptyword$,
and the singleton $\set{\emptyword}$ is the smallest nonempty
prefix-closed subset of $\Sigma^{*}$. An analogous observation applies
to suffix-closed sets, too.

We say that a subset $L$ of $\Sigma^{*}$ is \emph{subword-closed} if,
whenever $w_1, w_2, \ldots, w_n $ are elements of $\Sigma^{*}$, whose
product $w_1 w_2 \cdots w_n$ belongs to $L$, we have $w_i \in L$ for
all
$i = 1, \ldots, n$. This is equivalent to the condition that $L$ is
both prefix-closed and suffix-closed.

A \emph{cycle} in a binary graph $\Gamma = (q, \delta)$ is a path
$P = (q_0, w)$, such that $q_0 w = q_0$. We say that $\Gamma$ is
\emph{acyclic} if every cycle in $\Gamma$ is a trivial path.

\begin{proposition}
  \label{pro:PartialOrder}
  Let $\Gamma = (Q, \delta)$ be a binary graph. Then, the following
  statements are equivalent:
  \begin{enumerate}
  \item \label{item:Poset} The inheritance preorder $\leq$ on $Q$ (see
    Definition \ref{def:InheritancePreorder}) is a partial order.
  \item \label{item:PrefixClosed} For every node $q \in Q$, the
    isotropy submonoid $\Sigma^{*}_{q}$ is prefix-closed.
  \item \label{item:SubwordClosed} For every node $q \in Q$, the
    isotropy submonoid $\Sigma^{*}_{q}$ is subword-closed.
  \item \label{item:EmptyZeroOne} For every internal node $q \in Q$,
    the isotropy submonoid $\Sigma^{*}_{q}$ is one of the three
    submonoids $\set{\emptyword}$, $\set{0}^{*}$ and $\set{1}^{*}$.
  \item \label{item:Acyclic} The binary graph $\Gamma$ is acyclic.
  \end{enumerate}
\end{proposition}

\noindent \emph{Proof}. Suppose, first, that the inheritance preorder,
$\leq$, is anti-symmetric. Let $q \in Q$ be a node, let
$w \in \Sigma^{*}_{q}$, and let $v$ be a prefix of $w$, say $w = vz$,
for some word $z$. Then the node $q' = qv$ clearly satisfies the
inequality $q' \leq q$. On the other hand, $q = q'z$, so $q \leq q'$.
Therefore, by the anti-symmetry of $\leq$, we obtain $q' = q$. Since
$q = q' = qv$, the word $v$ belongs to $\Sigma^{*}_{q}$ and, hence,
that set is prefix-closed. We have, thus, proved that the condition
(\ref{item:Poset}) implies (\ref{item:PrefixClosed}).

Let us, now, suppose that $\Sigma^{*}_{q}$ is prefix-closed, and that
$w_1, \ldots, w_n$ are words, whose product $w = w_1 \cdots w_n$
belongs to $\Sigma^{*}_{q}$. Since the isotropy submonoid is
prefix-closed, it contains, for each $i = 1, \ldots, n$, the prefixes
$w_1 \cdots w_{i-1}$ and $w_1 \cdots w_i$, of $w$. Therefore,
\begin{displaymath}
  q w_i = (q w_1 \cdots w_{i-1}) w_i =
  q (w_1 \cdots w_{i-1} w_i) = q,
\end{displaymath}
i.e., $w_i \in \Sigma^{*}_{q}$, for all $i = 1, \ldots, n$. So,
$\Sigma^{*}_{q}$ is subword-closed. Conversely, every subword-closed
set is, obviously, prefix-closed. This proves that
(\ref{item:PrefixClosed}) and (\ref{item:SubwordClosed}) are
equivalent statements.

The conditions (\ref{item:SubwordClosed} and (\ref{item:EmptyZeroOne})
are equivalent, for the only subword-closed proper submonoids of
$\Sigma^{*}$ are $\set{\emptyword}$, $\set{0}^{*}$ and $\set{1}^{*}$.

Assume that $\Gamma$ satisfies the condition
(\ref{item:PrefixClosed}), and let $P = (q_0, w)$ be a cycle in
$\Gamma$.  Suppose that $P$ has positive length $n$, say,
$w = a_1 \cdots a_n$, where $a_i \in \Sigma$.
Let $q_i = q_0 a_1 \cdots a_n$, for $i = 1, \ldots, n$. Since $P$ is a
cycle, $q_0 w = q_0$, hence $w \in \Sigma^{*}_{q_0}$. The set
$\Sigma^{*}_{q_0}$ is prefix-closed, so it contains the prefix $a_1$
of $w$. Therefore, $q_1 = q_0 a_1 = q_0$, a contradiction since
$q_{i-1}$ is a parent of $q_i$ for all $i$. The cycle $P$, therefore,
has to be a trivial path, so (\ref{item:PrefixClosed}) implies
(\ref{item:Acyclic}).

Suppose, finally, that $\Gamma$ is an acyclic binary graph, and let
$q, q' \in Q$ be such that $q \leq q'$ and $q' \leq q$. Since
$q' \leq q$, by the definition of the partial order on $Q$, there
exists a path $P = (q_0, w)$ in $\Gamma$ such that
$q_0 = q$, and $q_0 w = q'$. Since $q \leq q'$, there exists a path
$P' = (q_0', w')$ such that $q_0' = q'$ and $q_0' w' = q$. The pair
$PP'$, defined by
$PP' = (q_0, ww')$, is a path, and is, in fact, a cycle in $\Gamma$.
Since $\Gamma$ is acyclic, the cycle $PP'$ must have length $0$, hence
$w = w' = \emptyword$. Therefore,
$q' = q_0 w = q_0 =q$ and, so, the preorder $\leq$ is anti-symmetric.
We have, now, shown that the condition (\ref{item:Acyclic}) implies
(\ref{item:Poset}). \qed

\section{Transition diagrams}
\label{sec:transitiondiagrams}

Transition diagrams are a device through which we connect our
description of binary graphs with their natural treatment in graph
theory. Let $\Gamma = (Q, \delta)$ be a binary graph. Then the
\emph{transition diagram} of $\Gamma$ is a labelled directed graph
$D (\Gamma)$, which is defined as follows:
\begin{enumerate}
\item The vertex set of $D (\Gamma)$ is defined to be $Q$.
\item If $(q, q', a) \in Q \times Q \times \Sigma$, if $q \neq q'$,
  and if $\delta (q, a) = q'$, then there is an edge labelled $a$ from
  $q$ to $q'$.
\end{enumerate}
Thus, there is an edge from $q$ to $q'$ if and only if $q$ is a parent
of $q'$. If $(q, a) \in Q \times \Sigma$, and if there is no edge from
$q$ labelled $a$, then it means that $\delta (q, a) = q$. In
particular, the vertices of $D (\Gamma)$ which have no edge starting
from them are precisely the leaves in $Q$.

\begin{example}
  \label{exa:trandiag1}
  The following directed graph is the transition diagram of the binary
  graph of Example \ref{exa:binarygraph1}.
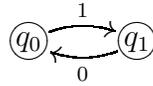
  \begin{displaymath}
    \entrymodifiers={+[o][F-]}
    \SelectTips{cm}{}
    \xymatrix{%
      q_0 \ar@/^/[r]^1 & q_1 \ar@/^/[l]^0
    }
  \end{displaymath}
\end{example}

\begin{example}
  \label{exa:trandiag2}  
  The transition diagram of the binary graph defined in Example
  \ref{exa:binarygraph2}, is as follows:
  \begin{displaymath}
    \SelectTips{cm}{}
    \xymatrix{%
      & & *+[o][F-]{q_0} \ar[dl]_0 \ar[dr]^1 & & \\
      & *+[o][F-]{q_1} \ar[dl]_0 \ar[dr]^1 & &
      *+[o][F-]{q_2} \ar[dl]_0 \ar[dr]^1 & \\
      *+[o][F-]{q_3} & & *+[o][F-]{q_4} & &
      *+[o][F-]{q_5} \ar[dl]_0 \\
      & & & *+[o][F-]{q_6} &
    }
  \end{displaymath}
\end{example}

\begin{example}
  \label{exa:trandiag3}
  The binary graph of Example \ref{exa:binarygraph3} has the following
  transition diagram:
  \begin{displaymath}
    \SelectTips{cm}{}
    
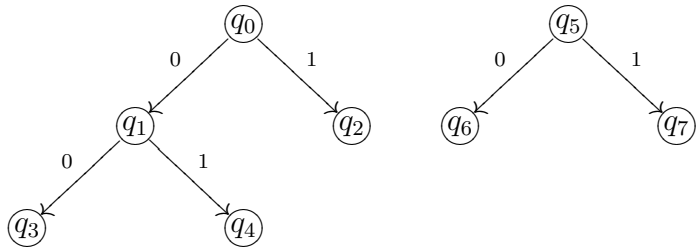
\xymatrix{%
      & & *+[o][F-]{q_0} \ar[dl]_0 \ar[dr]^1 & & & *+[o][F-]{q_5}
      \ar[dl]_0 \ar[dr]^1 & \\
      & *+[o][F-]{q_1} \ar[dl]_0 \ar[dr]^1 & & 
      *+[o][F-]{q_2} & *+[o][F-]{q_6} & & *+[o][F-]{q_7} \\
      *+[o][F-]{q_3} & & *+[o][F-]{q_4} & & & &
    }
  \end{displaymath}
\end{example}

\begin{example}
  Here is the transition diagram of the binary graph described in
  Example \ref{exa:binarygraph4}:
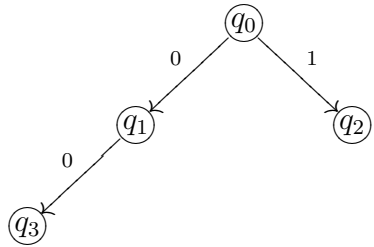
  \begin{displaymath}
    \SelectTips{cm}{}
    \xymatrix{%
      & & *+[o][F-]{q_0} \ar[dl]_0 \ar[dr]^1 & \\
      & *+[o][F-]{q_1} \ar[dl]_0 & & *+[o][F-]{q_2} \\
      *+[o][F-]{q_3} & & &
    }
\end{displaymath}
\end{example}

\section{Binary forests and trees}
\label{sec:binaryforests}

In this Section, we will define binary forests and trees, in our
formalism, and discuss their basic properties in terms of the
inheritance preorder and the action of bit strings.

\begin{definition}
  \label{def:BinaryForest}
  A \emph{binary forest} is a binary graph $\Gamma = (Q, \delta)$,
  which satisfies the following conditions:
  \begin{enumerate}
  \item \label{item:defPartialOrder} The inheritance preorder on $Q$
    is a partial order (see Definition \ref{def:InheritancePreorder}
    and Proposition \ref{pro:PartialOrder}).
  \item \label{item:SingleParent} If $q a = q' b$, where $a$ and $b$
    are letters in $\Sigma$, then either $q a = q$, or $q' b = q$, or
    $(q, a) = (q', b)$.
  \end{enumerate}
\end{definition}

\begin{example}
  \label{exa:BinaryForest}
  The condition (\ref{item:SingleParent}) in Definition
  \ref{def:BinaryForest} can be restated as the condition that every
  node can have at most one parent. For instance, the binary graph in
  Example \ref{exa:binarygraph3} satisfies this condition, and is a
  binary forest. The graph of Example \ref{exa:binarygraph2} does not
  satisfy the condition (\ref{item:SingleParent}). The graph in
  Example \ref{exa:binarygraph1} does not satisfy the condition
  (\ref{item:defPartialOrder}). Note that every subgraph of a binary
  forest is, again, a binary forest. We will sometimes refer to such a
  subgraph as a \emph{subforest}.
\end{example}

We say that a path $P = (q_0, w)$ is a \emph{suffix} of a path
$P' = (q_0', w')$, if there exists a word $v \in \Sigma^{*}$, such
that $w' = v w$, and $q_0 = q_0' v$. Note that if this is the case,
then $P$ and $P'$ end at the same node.

\begin{proposition}
  \label{pro:SuffixPath}
  A binary graph $\Gamma = (Q, \delta)$ satisfies the condition
  (\ref{item:SingleParent}) of Definition \ref{def:BinaryForest}, if
  and only if it satisfies the following condition:
  \begin{enumerate}
  \item[($2'$)] If $P$ and $P'$ are paths in $\Gamma$, which end at
    the same node, then either $P$ is a suffix of $P'$, or $P'$ is a
    suffix of $P$.
  \end{enumerate}
\end{proposition}

\noindent \emph{Proof}. Suppose that $\Gamma$ satisfies the condition
($2'$) of the Proposition. Let $q, a, q', b$ be as in
(\ref{item:SingleParent}) of Definition \ref{def:BinaryForest}, and
suppose that $q a \neq q$ and $q' b \neq q'$. Then $P = (q, a)$ and
$P' = (q', b)$ are paths in $\Gamma$ which end at the same node.
Therefore, by ($2'$), one of them is a suffix of the other. As the
both have the same length, $1$, this means that they are equal, i.e.,
$q = q'$ and $a = b$, so, ($2'$) implies condition
(\ref{item:SingleParent}) of Definition \ref{def:BinaryForest}.

Conversely, assume that $\Gamma$ satisfies the condition
(\ref{item:SingleParent}) of Definition \ref{def:BinaryForest}. Let
$P = (q_0, w)$ and $P' = (q_0', w')$ be two paths, which end at the
same node, i.e., $q_0 w = q_0' w'$. It is obvious that if $w$ is the
empty word, then $P$ is a suffix of $P'$ and, similarly, if $w'$ is
the empty word, then $P'$ is a suffix of $P$. So, let us assume that
both $w$ and $w'$ have positive length, say, $w = a_1 \cdots a_m$ and
$w' = a_1' \cdots a_n'$. Let $q_i = q_0 a_1 \cdots a_i$ and
$q_j = q_0' a_1' \cdots a_j'$, for $i = 1, \ldots, m$ and
$j = 1, \ldots, n$. Suppose, without loss of generality, that
$m \leq n$. Since $P$ and $P'$ are paths, we have
$q_{m-1} \neq q_m$ and $q_{n-1}' \neq q_n'$. So, applying the
condition (\ref{item:SingleParent}) to the equation
$q_{m-1} a_m = q_{n-1}' a_n'$, we get that
$q_{m-1} = q_{n-1}'$ and $a_m = a_n$. Considering, now, the paths
$(q_0, a_1 \cdots a_{m-1})$ and $(q_0', a_1' \cdots a_{n-1}')$, and
using induction, we see that $q_i = q_{n-m+i}'$ for all
$i = 0, 1, \ldots, m$, and $a_i = a_{n-m+i}'$ for all
$i = 1, \ldots, m$. Setting $v = a_1' \cdots a_{n-m}'$, we have
$w' = v w$, and $q_0' v = q_{n-m}' = q_0$. Therefore, $P$ is a
suffix of $P'$. \qed

\begin{corollary}
  \label{cor:UniquePath}
  Let $\Gamma = (Q, \delta)$ be a binary forest. Then, for every pair
  of nodes $q, q' \in Q$ such that $q' \leq q$, there is a unique path
  from $q$ to $q'$.
\end{corollary}

\noindent \emph{Proof}. Let $P = (q, w)$ and $P' = (q, w')$ be two
paths from $q$ to $q'$.  Then by Proposition \ref{pro:SuffixPath}, one
of these two paths is a suffix of the other, say, $P$ is a suffix of
$P'$. Thus, there exists a word $v \in \Sigma^{*}$, such that
$w' = v w$, and $q = q v$. Since $P'$ is a path, so is the pair
$Q = (q, v)$, and clearly $Q$ starts and ends at $q$. By Proposition
\ref{pro:PartialOrder}, $\Gamma$ is acyclic, hence the cycle $Q$ is a
trivial path. Therefore, $v = \emptyword$, which implies that $w' = w$
and $P' = P$. \qed

Let $X$ be a set with a preorder $\leq$, and let $S$ be a subset of
$X$. Recall that an element $x_0$ of $S$ is said to be a
\emph{maximum} element of $S$ if $x \leq x_0$ for all $x \in S$, and
that $x_0$ is said to be a \emph{maximal} element of $S$ if, whenever
$x \in S$ and $x_0 \leq x$, we have $x = x_0$.

\begin{remark}
  \label{rem:UniqueMaximum}
  Note that if the preorder $\leq$ on $X$ is a partial order, then a
  subset of $X$ has at most one maximum. Also, in any preorder, if a
  set has a maximal element and a maximum element, then they must
  coincide, and the set has a unique maximum.
\end{remark}

\begin{proposition}
  \label{pro:MaximalOrbits}
  Consider a binary graph $\Gamma = (Q, \delta)$, and the inheritance
  preorder, $\leq$, on $Q$ (see Definition
  \ref{def:InheritancePreorder}). Then, the following statements are
  true:
  \begin{enumerate}
  \item \label{item:MaximalOrbits} The set of nodes, $Q$, has a
    maximal element, and is the union of the $\Sigma^{*}$-orbits of
    all the maximal elements.
  \item \label{item:ForestMaximalOrbits} If $\Gamma$ is a binary
    forest, and if $q_1$ and $q_2$ are distinct maximal elements of
    $Q$, then their orbits,
    $\Sigma^{*} (q_1)$ and $\Sigma^{*} (q_2)$ are disjoint.
  \end{enumerate}
\end{proposition}

\noindent \emph{Proof}. Let $q_0 \in Q$ be a node. For each
$i \in \N$, inductively define a node $q_i$, as follows:
$q_i = q_{i-1}$ if $q_{i-1}$ is maximal, and if $q_{i-1}$ is not
maximal, then $q_i$ is some choice of a node $q$ such that
$q_{i-1} < q$. Since $Q$ is a finite set, the sequence
$\set{q_i}_{i=0}^{\infty}$ is eventually constant, i.e., there exists
$n \geq 0$ such that $q_i = q_n$ for all $i \geq n$. Then, $q_n$ is a
maximal node, and $q_0 \leq q_n$, i.e.,
$q_0 \in \Sigma^{*} (q_n)$. This proves that
(\ref{item:MaximalOrbits}) is true for any binary graph.

Suppose, next, that $q_1$ and $q_2$ are maximal elements in a binary
forest $\Gamma$, and that
$q \in \Sigma^{*} (q_1) \cap \Sigma^{*} (q_2)$. Then, there exist a
path $P_1$ from $q_1$ to $q$, and a path $P_2$ from $q_2$ to $q$.
Since $P_1$ and $P_2$ end at the same node, one of them is a suffix of
the other, say, $P_1$ is a suffix of $P_2$. It follows that
$q_1 \in \Sigma^{*} (q_2)$, so, $q_1 \leq q_2$. By maximality of
$q_1$, we obtain that $q_1 = q_2$. Therefore, the orbits of two
distinct maximal elements in $\Gamma$ cannot intersect. \qed

\begin{definition}
  \label{def:Connected}
  Let $\Gamma = (Q, \delta)$ be a binary graph. We say that a subset
  $S$ of $Q$ is a \emph{connected set} if it has a maximum element. A
  \emph{connected component} of $\Gamma$ is a maximal connected subset
  of $Q$. We say that $\Gamma$ is \emph{connected} if its set of
  nodes, $Q$, is a connected set.
\end{definition}

If $\Gamma = (Q, \delta)$ is a connected binary graph, then $Q$ has a
unique maximum. For, by (\ref{item:MaximalOrbits}) of Proposition
\ref{pro:MaximalOrbits}, $Q$ has a maximal element. By Remark
\ref{rem:UniqueMaximum}, this element has to be the unique maximum of
$Q$. We call it the \emph{root} of $\Gamma$.

Suppose that $S$ is a connected component of a binary graph
$\Gamma = (Q, \delta)$, and let $q_0$ be a maximum element of $S$. By
(\ref{item:MaximalOrbits}) of Proposition \ref{pro:MaximalOrbits},
there exists a maximal node $q_1$ in $Q$, such that $q_0 \leq q_1$.
For every element $q \in S$, we have $q \leq q_0 \leq q_1$, hence $S$
is contained in the orbit $\Sigma^{*} (q_1)$. This orbit is a
connected set, having $q_1$ as a maximum. Therefore, by the maximality
of $S$, we have $S = \Sigma^{*} (q_1)$. Since $S$ has a maximum
element $q_0$ and a maximal element $q_1$, by Remark
\ref{rem:UniqueMaximum}, we have $q_0 = q_1$, and $q_0$ is the unique
maximum of $S$. In other words, every connected component $S$ of a
binary graph has a unique maximum $q_0$, and this maximum is a maximal
element in the full set of nodes; moreover, $S = \Sigma^{*} (q_0)$.
Therefore, by Remark \ref{rem:OrbitGraph}, we have a binary subgraph
$\Gamma_{q_0} = (\Sigma^{*} (q_0), \delta_{q_0})$ of $\Gamma$, having
$S = \Sigma^{*} (q_0)$ as its set of nodes. We denote the graph
$\Gamma_{q_0}$ by $\Gamma_S$ and, by abuse of language, we refer to
$\Gamma_S$ also as a connected component of $\Gamma$.

\begin{proposition}
  \label{pro:ConnectedComponents}
  Let $\Gamma = (Q, \delta)$ be a binary graph. Then, the orbit of any
  maximal node $q_0 \in Q$ is a connected component of $\Gamma$.
  Conversely, every connected component $S$ of $\Gamma$ has a unique
  maximum $q_0$; this maximum is a maximal element in the full set of
  nodes, and $S$ is the orbit of $q_0$.
\end{proposition}

\noindent \emph{Proof}. Suppose that $q_0$ is a maximal element of
$Q$. Obviously, $q_0$ is a maximum element of $\Sigma^{*} (q_0)$, so
the orbit of $q_0$ is a connected set. Let $S$ be a connected subset
of $Q$, which contains $\Sigma^{*} (q_0)$, and denote by $q_1$ a
maximum element of $S$.  Then, $q_0 \leq q_1$, so, by maximality of
$q_0$, we have $q_0 = q_1$.  Since $q_0$ is a maximum of $S$, we have
$q \leq q_0$ for all $q \in S$, i.e.,
$q \in \Sigma^{*} (q_0)$ for all $q \in S$. Therefore,
$S = \Sigma^{*} (q_0)$, hence, $\Sigma^{*} (q_0)$ is a maximal
connected subset of $Q$. This proves the first assertion of the
Proposition. We have already shown, above, the truth of the second
assertion. \qed

\begin{remark}
  \label{rem:LeftRightSubgraphs}
  If $\Gamma = (Q, \delta)$ is a binary graph, then, for any node
  $q_0 \in Q$, we define
  \begin{align*}
    Q_{\lambda} (q_0) & = \set{ q_0 0 v \suchthat v \in \Sigma^{*}},
    \quad \text{and} \\
    Q_{\rho} (q_0) & = \set{ q_0 1 v \suchthat v \in \Sigma^{*}}.
  \end{align*}
  Since these two sets are $\Sigma^{*}$-invariant, they define, by
  Remark \ref{rem:morphismbinarygraph}, two subgraphs
  $\Gamma_{\lambda} (q_0) = (Q_{\lambda} (q_0), \delta)$ and
  $\Gamma_{\rho} (q_0) = (Q_{\rho} (q_0), \delta)$, where, by abuse of
  notation, we denote the restriction of $\delta$ to a subset of
  $Q \times \Sigma$ by the same symbol. It is clear that
  $\delta (q_0, 0)$ is a maximum element of $Q_{\lambda} (q_0)$, so
  $\Gamma_{\lambda} (q_0)$ is a connected graph. We call it the
  \emph{left subgraph} at $q_0$.  Similarly, $\Gamma_{\rho} (q_0)$ is
  also a connected graph, having $\delta (q_0, 1)$ as the maximum
  element. We call it the
  \emph{right subgraph} at $q_0$.
\end{remark}

\begin{definition}
  \label{def:BinaryTree}
  A connected binary forest is called a \emph{binary tree}. In this
  case, the set $Q$ has a unique maximum (see Remark
  \ref{rem:UniqueMaximum}), and this element is called the \emph{root}
  of $\Gamma$.
\end{definition}

\begin{example}
  \label{exa:BinaryTree}
  The binary graph in Example \ref{exa:binarygraph3} is not connected,
  because the subset $\set{q_0, q_5}$ of $Q$ does not have an upper
  bound. In this graph, $q_0$ and $q_5$ are maximal nodes. In Example
  \ref{exa:binarygraph1}, there are no maximal nodes, but both the
  nodes $q_0$ and $q_1$ are maximum nodes. The graph in Example
  \ref{exa:binarygraph4} is a binary tree, with root $q_0$.
\end{example}

Let $\Gamma = (Q, \delta)$ is a binary forest. Then, every connected
component of $\Gamma$ is a binary tree. As we had observed in Example
\ref{exa:BinaryForest}, every subgraph of $\Gamma$ is also a binary
forest. In particular, for any node $q_0 \in Q$, the left and right
subgraphs at $q_0$ (See Remark \ref{rem:LeftRightSubgraphs}), namely
$\Gamma_{\lambda} (q_0)$ and $\Gamma_{\rho} (q_0)$, are binary
forests. Since they are, always, connected graphs, it follows that
$\Gamma_{\lambda} (q_0)$ and $\Gamma_{\rho} (q_0)$ are, in fact,
binary trees. We call them, respectively, the \emph{left subtree} and
the \emph{right subtree} at $q_0$.

\begin{proposition}
  \label{pro:DisjointSubtrees}
  Let $\Gamma = (Q, \delta)$ be a binary forest, and suppose that
  $q_0$ is a node in $Q$ which has trivial isotropy. Then, the left
  and right subtrees of $q_0$ are disjoint, i.e.,
  $Q_{\lambda} (q_0) \cap Q_{\rho} (q_0) = \emptyset$.
\end{proposition}

\noindent \emph{Proof}. Suppose that
$Q_{\lambda} (q_0) \cap Q_{\rho} (q_0) \neq \emptyset$, say
$q_0 0 v = q_0 1 w$ for some words $v, w \in \Sigma^{*}$. Let
$q' = q_0 0 v$. Then, $q'$ is a descendant of both $q_0 0$ and
$q_0 1$. Therefore, by Proposition \ref{pro:ancestor}, there exist
paths $P_0 = (q_0 0, z_0)$ and $P_1 = (q_0 1, z_1)$ from $q_0 0$ and
$q_0 1$, respectively, to $q'$. Since $q_0$ has trivial isotropy,
$q_0 0 \neq q_0$ and $q_0 1 \neq q_0$. Therefore, $Q_0 = (q_0, 0z_0)$
and $Q_1 = (q_0, 1z_1)$ are also paths. Since both of them start at
the same point $q_0$, and end at the same point $q'$, by Corollary
\ref{cor:UniquePath}, $Q_0 = Q_1$, i.e., $0z_0 = 1z_1$. Since this is
an impossibility, we conclude that our hypothesis, that the left and
right subtrees of $q_0$ overlap, has to be false. \qed

\section{Binary automata and transition systems}
\label{sec:binaryautomata}

We will, now, recall the notion of binary automata, and describe a
canonical functor from the category of these automata, to that of
binary graphs with a distinguished node. We will also relate binary
graphs to transition systems.

A \emph{binary automaton} is a finite automaton
$M = (Q, \Sigma, \delta, q_0, F)$, where the input alphabet $\Sigma$
equals the set $\set{0,1}$. (We follow the notation of
\cite[Chapter 1]{Hopcroft:1979:IAT} for finite automata.)  Let us
define a state $q$ to be \emph{stationary} if $\delta (q, a) = q$ for
all $a \in \Sigma$.  We will say that a binary automaton is
\emph{admissible}, in case the set $F$ of final states coincides with
the set of stationary states.

Let $M = (Q, \Sigma, \delta, q_0, F)$ and
$M' = (Q', \Sigma, \delta', q_0', F')$ be binary automata. A
\emph{morphism} from $M$ to $M'$ is a function
$f : Q \rightarrow Q'$, such that the following conditions hold:
\begin{enumerate}
\item $f (\delta (q, a)) = \delta' (f (q), a)$ for all $q \in Q$ and
  $a \in \Sigma$;
\item $f (q_0) = q_0'$; and
\item $f (F) \subset F'$.
\end{enumerate}
We say that $M$ is a \emph{binary subautomaton} of $M'$, if $Q \subset
Q'$, and if the inclusion $i : Q \hookrightarrow Q'$ is a morphism of
binary automata.  This is equivalent to the assertion that $Q \subset
Q'$, that $\delta$ is the restriction of $\delta'$, that $q_0' \in Q$,
and that
$F \subset F'$.

A \emph{pointed binary graph} is a pair $P = (\Gamma, q_0)$, where
$\Gamma = (Q, \delta)$ is a binary graph, and $q_0 \in Q$ is a
distinguished node, called the \emph{base node}. Let
$P = (\Gamma, q_0)$ and $P' = (\Gamma', q_0')$ be pointed binary
graphs. A \emph{morphism} from $P$ to $P'$ is a morphism of binary
graphs $f : \Gamma \rightarrow \Gamma'$, such that $f (q_0) = q_0'$.
We say that $P$ is a \emph{pointed binary subgraph} of $P'$ if
$\Gamma$ is a binary subgraph of $\Gamma'$, and if the canonical
monomorphism of binary graphs, $i : \Gamma \hookrightarrow \Gamma'$,
is a morphism of pointed binary graphs. If $\Gamma = (Q, \delta)$ and
$\Gamma' = (Q', \delta')$, this is equivalent to saying that
$Q \subset Q'$, that $\delta$ is the restriction of $\delta'$, and
that $q_0' \in Q$.

Given a binary automaton $M = (Q, \Sigma, \delta, q_0, F)$, we get a
pointed binary graph $P_M = (\Gamma_M, q_0)$, where
$\Gamma_M = (Q, \delta)$. The leaves of $\Gamma_M$ are precisely the
stationary states of $M$. In particular, if the automaton $M$ is
admissible, the leaves of $\Gamma_M$ are the same as the final states
of $M$. If $M = (Q, \Sigma, \delta, q_0, F)$, and
$M' = (Q', \Sigma, \delta', q_0', F')$ are two binary automata, and if
$f : M \rightarrow M'$ is a morphism, then the underlying function $f
: Q \rightarrow Q'$ defines a morphism of pointed binary graphs
$P_f : P_M \rightarrow P_{M'}$. We thus obtain a functor $p$ from the
category $\mathcal{BA}$ of binary automata to the category
$\mathcal{BG}_{*}$ of pointed binary graphs
\begin{gather*}
  p : \Ob{\mathcal{BA}} \rightarrow \Ob{\mathcal{BG}_{*}}, \quad 
  p (M) = P_M, \quad \text{and} \\
  p : \Hom{M}{M'} \rightarrow \Hom{P_M}{P_{M'}}, \quad p (f) = P_f.
\end{gather*}
We will see, in Section \ref{sec:binaryfibredcat}, that this morphism
has quite an interesting structure.

\begin{remark}
  \label{rem:admissible}
  The leaves of $\Gamma_M$ are precisely the stationary states of $M$.
  In particular, if the automaton $M$ is admissible, the leaves of
  $\Gamma_M$ are the same as the final states of $M$. Moreover, the
  language $L (M) \subset \Sigma^{*}$ of an admissible automaton $M$
  is a monoidal ideal, i.e., if $w \in L (M)$ and $z \in \Sigma^{*}$,
  then $w z \in L (M)$.
\end{remark}

In the theory of concurrent computation, one has the notion of a
transition system (see \cite[Chapter 2]{Winskel:1993:MC}). A
\emph{transition system} is a datum
$T = (S, i, L, \mathit{Tran})$, where
\begin{enumerate}
\item $S$ is a set, whose elements are called \emph{states}.
\item $i$ is a distinguished element of $S$, called the
  \emph{initial state}.
\item $L$ is a set, whose elements are called \emph{labels}.
\item $\mathit{Tran}$ is a subset of $S \times L \times S$, called the
  \emph{transition relation}.
\end{enumerate}

Let $T = (S, i, L, \mathit{Tran})$ and
$T' = (S', i', L', \mathit{Tran}')$ be transition systems. A
\emph{morphism} from $T$ to $T'$ is a pair $(\sigma, \lambda)$, where
$\sigma : S \rightarrow S'$ is a function and $\lambda : L
\dashrightarrow L'$ is a partial function (i.e., a function from a
subset of $L$ to $L'$), satisfying the following conditions:
\begin{enumerate}
\item $\sigma (i) = i'$.
\item If $(s_1, a, s_2) \in \mathit{Tran}$ and if $\lambda (a)$ is
  defined, then
  $(\sigma (s_1), \lambda (a), \sigma (s_2)) \in \mathit{Tran}'$.
\item If $(s_1, a, s_2) \in \mathit{Tran}$ and if $\lambda (a)$ is not
  defined, then $\sigma (s_1) = \sigma (s_2)$.
\end{enumerate}

If $P = (\Gamma, q_0)$ is a pointed binary graph, with
$\Gamma = (Q, \delta)$, we obtain a transition system
$T_P = (Q, q_0, \Sigma, \mathit{Gr} (\delta))$, where
$\mathit{Gr} (\delta) \subset Q \times \Sigma \times Q$ denotes the
graph of $\delta$, i.e.,
\begin{displaymath}
  \mathit{Gr} (\delta) = \set{%
    (q_1, a, q_2) \in Q \times \Sigma \times Q \suchthat
    q_2 = \delta (q_1, a)}.
\end{displaymath}
Moreover, if $f : P \rightarrow P'$ is a morphism from a pointed
binary  graph $P = (\Gamma, q_0)$ to a pointed binary graph
$P' = (\Gamma', q_0')$, then the pair $(f, \id{\Sigma})$ is a morphism
of the associated transition systems, $T_P \rightarrow T_{P'}$. This
assignment provides us a functor
$\mathcal{BG}_{*} \rightarrow \mathcal{TS}$, from the category of
pointed binary graphs, to the category of transition systems. This
functor is an isomorphism from $\mathcal{BG}_{*}$ to a subcategory of
$\mathcal{TS}$. However, $\mathcal{BG}_{*}$ is not a full subcategory
of $\mathcal{TS}$, as can be seen using simple examples of pointed
binary graphs with two nodes.

\section{Fibred categories}
\label{sec:fibredcategories}

The original source for fibred categories is Grothendieck's article in
SGA~1 \cite[Expos\'{e}~VI]{Grothendieck:1971:REG}. We also refer the
reader to the notes of Streicher \cite{Streicher:1999:FCJ}. We recall
below the basic notions about fibred categories. To aid the reader,
and to make this article self-contained, we provide in this Section
and the following two Sections, proofs of all propositions we need
regarding fibred categories. The material is, perhaps, standard for
specialists in certain areas of algebraic geometry and category
theory. Such readers may proceed directly to Section
\ref{sec:binaryfibredcat}.

Let $\mathcal{E}$ be a category. An $\mathcal{E}$-\emph{category}, or
a \emph{category over} $\mathcal{E}$, is a pair $(\mathcal{F}, p)$,
where $\mathcal{F}$ is a category, and
$p : \mathcal{F} \rightarrow \mathcal{E}$ is a functor. We usually
drop the functor $p$ from the notation, and say that $\mathcal{F}$ is
an $\mathcal{E}$-category. Sometimes, we also say that $p :
\mathcal{F} \rightarrow \mathcal{E}$ is an $\mathcal{E}$-category. An
$\mathcal{E}$-\emph{functor} from an $\mathcal{E}$-category
$(\mathcal{F}, p)$ to an $\mathcal{E}$-category $(\mathcal{G}, q)$ is
a functor
$f : \mathcal{F} \rightarrow \mathcal{G}$, such that $q \circ f = p$.
In particular, a \emph{section} of $(\mathcal{F}, p)$ is a functor
$s : \mathcal{E} \rightarrow \mathcal{F}$ such that
$p \circ s = \id{\mathcal{E}}$.

Let $\mathcal{F}$ be an $\mathcal{E}$-category, and let
$S \in \Ob{\mathcal{E}}$. The \emph{categorical fibre} of
$\mathcal{F}$ at $S$ is the subcategory $\mathcal{F}_S$ of
$\mathcal{F}$ defined as follows:
\begin{enumerate}
\item $\Ob{\mathcal{F}_S}$ is the collection of all
  $\xi \in \Ob{\mathcal{F}}$ such that $p (\xi) = S$.
\item If $\xi, \xi' \in \Ob{\mathcal{F}_S}$, a morphism from $\xi$ to
  $\xi'$ in $\mathcal{F}_S$ is a morphism
  $u : \xi \rightarrow \xi'$ in $\mathcal{F}$, such that
  $p (u) = \id{S}$. We will call such a morphism $u$ an
  $S$-\emph{morphism}, and we denote by $\Hom[S]{\xi}{\xi'}$, the set
  of all $S$-morphisms from $\xi$ to $\xi'$.
\end{enumerate}

We generalize the above notion of an $\mathcal{F}_S$-morphism as
follows. Suppose $f : T \rightarrow S$ is a morphism in $\mathcal{E}$,
let $\eta \in \Ob{\mathcal{F}_T}$, and let
$\xi \in \Ob{\mathcal{F}_S}$. An $f$-\emph{morphism} from $\eta$ to
$\xi$ is a morphism $u : \eta \rightarrow \xi$ in $\mathcal{F}$, such
that $p (u) = f$. We denote the collection of all $f$-morphisms from
$\eta$ to $\xi$ by $\Hom[f]{\eta}{\xi}$.

\begin{definition}
  \label{def:cartesianmorphism}
  (See
  \cite[Expos\'{e} VI, D\'{e}finition 5.1]{Grothendieck:1971:REG}.)
  Let $(\mathcal{F}, p)$ be an $\mathcal{E}$-category. Let
  $\alpha : \eta \rightarrow \xi$ be a morphism in $\mathcal{F}$, and
  let $S = p (\xi)$, $T = p (\eta)$ and $f = p (\alpha)$. We say that
  $\alpha$ is a \emph{Cartesian} morphism if for every object
  $\eta' \in \Ob{\mathcal{F}_T}$, and for every $f$-morphism
  $u : \eta' \rightarrow \xi$, there exists a unique $T$-morphism
  $\overline{u} : \eta' \rightarrow \eta$ such that
  $u = \alpha \circ \overline{u}$. The pair $(\eta, \alpha)$ is called
  an \emph{inverse image} of $\xi$ by $f$.
\end{definition}

The situation in the above definition can be described as in the
following diagram:
\begin{displaymath}
  \xymatrix{%
    \eta' \ar[drrr]^u \ar@{-->}[dr]_{\overline{u}} & & & \\
    & \eta \ar[rr]_{\alpha} & & \xi \\
    & & \ar@{-}[d]^{p} & \\
    T \ar[drrr]^f \ar[dr]_{\id{T}} & & & \\
    & T \ar[rr]_f & & S
  }
\end{displaymath}
Thus, a morphism $\alpha : \eta \rightarrow \xi$ in $\mathcal{F}$ is
Cartesian if and only if for every $f$-morphism $u : \eta' \rightarrow
\xi$, the factorization $f = f \circ \id{T}$ of $f = p (\alpha)$ in
$\mathcal{E}$ lifts uniquely to a factorization
$u = \alpha \circ \overline{u}$ of $u$ in $\mathcal{F}$. In other
words, $\alpha : \eta \rightarrow \xi$ is Cartesian if and only if for
every $\eta' \in \Ob{\mathcal{F}_T}$, the function
\begin{align*}
  \alpha_{*} : \Hom[T]{\eta'}{\eta} &
  \rightarrow \Hom[f]{\eta'}{\xi} \\
  v & \mapsto \alpha \circ v
\end{align*}
is a bijection. We will use the notation $\alpha_{*}$ often, and
record it for clarity.

\begin{notation}
  Let $\mathcal{F}$ be a category, and let
  $\alpha : \eta \rightarrow \xi$ be a morphism in $\mathcal{F}$.
  Then, for every $\zeta \in \Ob{\mathcal{F}}$, we denote by
  $\alpha_{*}$ the function
  \begin{align*}
    \Hom{\zeta}{\eta} & \rightarrow \Hom{\zeta}{\xi}, \\
    v & \mapsto \alpha \circ v.
  \end{align*}
  We denote the restriction of $\alpha_{*}$ to any subset $F$ of
  $\Hom{\zeta}{\eta}$, also, by $\alpha_{*}$.
\end{notation}

Let $f : T \rightarrow S$ be a morphism in $\mathcal{E}$, and let
$\xi \in \Ob{\mathcal{F}_S}$. Suppose that $(\eta, \alpha)$ and
$(\eta', \alpha')$ are two inverse images of $\xi$ by $f$. Then, by
definition, there exists a unique $T$-isomorphism
$\overline{\alpha'} : \eta' \rightarrow \eta$ such that the diagram
\begin{displaymath}
  \xymatrix{%
    \eta' \ar[drrr]^{\alpha'}
    \ar@{-->}[dr]_{\overline{\alpha'}} & & & \\
    & \eta \ar[rr]_{\alpha} & & \xi \\
  }
\end{displaymath}
commutes. Thus, an inverse image of $\xi$ by $f$, if it exists, is
unique up to a canonical $T$-isomorphism. Suppose that an inverse
image of $\xi$ by $f$ exists, and that we have made a choice of such
an inverse image. In such a situation, we often denote the chosen
inverse image by $(f_{\mathcal{F}}^{*} \xi, \alpha_f (\xi))$, or
simply by
$(f^{*} \xi, \alpha_f (\xi))$. Further, by abuse of language, we then
call $(f^{*} \xi, \alpha_f (\xi))$ \emph{the} inverse image of $\xi$
by $f$. If an inverse image of $\xi$ by $f$ exists for all morphisms
$f: T \rightarrow S$ in $\mathcal{E}$, and for all
$\xi \in \Ob{\mathcal{F}_S}$, then we say that
\emph{the inverse image functor by $f$ in $\mathcal{F}$ exists}.  If
this is indeed the case, and if we have chosen, for all $f$ and for
all $\xi$, an inverse image of $\xi$ by $f$, then, the assignment
\begin{gather*}
  f^{*} : \Ob{\mathcal{F}_S} \rightarrow \Ob{\mathcal{F}_T}, \quad
  \xi \mapsto f^{*} \xi \\
  f^{*} : \Hom[S]{\xi}{\xi'} \rightarrow \Hom[T]{f^{*}\xi}{f^{*}\xi'},
  \quad u \mapsto f^{*} u,
\end{gather*}
is a functor from $\mathcal{F}_S$ to $\mathcal{F}_T$. Here,
$f^{*} u : f^{*} \xi \rightarrow f^{*} \xi'$ is the unique
$T$-morphism induced by the $f$-morphism
$u \circ \alpha_f (\xi) : f^{*} \xi \rightarrow \xi'$, using the
universal property of $f^{*} \xi'$.

\begin{definition}
  \label{def:fibredcategory}
  (See
  \cite[Expos\'{e} VI, D\'{e}finition 6.1]{Grothendieck:1971:REG}.)
  Let $\mathcal{E}$ be a category. We say that an
  $\mathcal{E}$-category $(\mathcal{F}, p)$ is a
  \emph{fibred category} over $\mathcal{E}$ if it satisfies the
  following two conditions:
  \begin{enumerate}
  \item[$\mathbf{Fib_I}$] For every morphism
    $f : T \rightarrow S$ in $\mathcal{E}$, the inverse image functor
    by $f$ in $\mathcal{F}$ exists.
  \item[$\mathbf{Fib_{II}}$] If $\alpha: \eta \rightarrow \xi$ and
    $\beta : \theta \rightarrow \eta$ are Cartesian morphisms in
    $\mathcal{F}$, then their composition,
    $\alpha \circ \beta : \theta \rightarrow \xi$, is also a Cartesian
    morphism.
  \end{enumerate}
\end{definition}

We will see presently, in Proposition \ref{pro:stronglycartfibred},
that in fibred categories, Cartesian morphisms satisfy a certain
stronger condition.

\section{A criterion for fibredness}
\label{sec:fibredcatcriterion}

We will now present a criterion for a category $\mathcal{F}$ over
$\mathcal{E}$ to be fibred. This criterion is sometimes taken to be a
definition of fibred categories, as in
\cite[Definitions 2.1 and 2.2]{Streicher:1999:FCJ}.

\begin{definition}
  \label{def:stronglycartesian}
  (See \cite[Appendix B, page 163]{Winskel:1993:MC}.) Let
  $(\mathcal{F}, p)$ be an $\mathcal{E}$-category. Let
  $\alpha : \eta \rightarrow \xi$ be a morphism in $\mathcal{F}$, and
  let $S = p (\xi)$, $T = p (\eta)$ and $f = p (\alpha)$. We say that
  $\alpha$ is a \emph{strongly Cartesian} morphism if for every
  morphism $g : U \rightarrow T$ in $\mathcal{E}$, for every object
  $\zeta \in \Ob{\mathcal{F}_U}$, and for every $(f \circ g)$-morphism
  $u : \zeta \rightarrow \xi$, there exists a unique $g$-morphism
  $\overline{u} : \zeta \rightarrow \eta$ such that
  $u = \alpha \circ \overline{u}$.
\end{definition}

The situation in the above definition can be described as in the
following diagram:
\begin{displaymath}
  \xymatrix{%
    \zeta \ar[drrr]^u \ar@{-->}[dr]_{\overline{u}} & & & \\
    & \eta \ar[rr]_{\alpha} & & \xi \\
    & & \ar@{-}[d]^{p} & \\
    U \ar[drrr]^{f \circ g} \ar[dr]_{g} & & & \\
    & T \ar[rr]_f & & S
  }
\end{displaymath}
Thus, a morphism $\alpha : \eta \rightarrow \xi$ in $\mathcal{F}$ is
strongly Cartesian if and only if for every morphism
$g : U \rightarrow T$ in $\mathcal{E}$, for every
$\zeta \in \Ob{\mathcal{F}_U}$, and for every $(f \circ g)$-morphism
$u : \zeta \rightarrow \xi$, the factorization $f \circ g = f \circ g$
of $f \circ g$ in $\mathcal{E}$ lifts uniquely to a factorization
$u = \alpha \circ \overline{u}$ of $u$ in $\mathcal{F}$. In other
words, $\alpha : \eta \rightarrow \xi$ is strongly Cartesian if and
only if for every morphism $g : U \rightarrow T$ in $\mathcal{E}$, and
for every $\zeta \in \Ob{\mathcal{F}_U}$, the function
\begin{align*}
  \alpha_{*} : \Hom[g]{\zeta}{\eta} &
  \rightarrow \Hom[f \circ g]{\zeta}{\xi} \\
  v & \mapsto \alpha \circ v
\end{align*}
is a bijection.

Strongly Cartesian morphisms behave well under composition.

\begin{proposition}
  \label{pro:stronglycartcompose}
  (See \cite[Appendix B, page 163]{Winskel:1993:MC}.) Let
  $\mathcal{E}$ be a category, and $(\mathcal{F}, p)$ an
  $\mathcal{E}$-category. Suppose $\alpha : \eta \rightarrow \xi$ and
  $\beta : \zeta \rightarrow \eta$ are strongly Cartesian morphisms in
  $\mathcal{F}$. Then, their composition
  $\alpha \circ \beta : \zeta \rightarrow \xi$ is also strongly
  Cartesian.
\end{proposition}

\noindent \emph{Proof.} Let $f : T \rightarrow S$ and
$g : U \rightarrow T$ denote $p (\alpha)$ and $p (\beta)$,
respectively. Let $h : V \rightarrow U$ be a morphism in
$\mathcal{E}$, and let $\theta \in \Ob{\mathcal{F}_V}$. Then, the
diagram
\begin{displaymath}
  \xymatrix{
    \Hom[h]{\theta}{\zeta} \ar[rrr]^{\beta_{*}}
    \ar[ddrrr]_{(\alpha \circ \beta)_{*}} & & &
    \Hom[g \circ h]{\theta}{\eta} \ar[dd]^{\alpha_{*}} \\
    & & & \\
    & & & \Hom[f \circ g \circ h]{\theta}{\xi}
  }
\end{displaymath}
commutes. Since $\alpha$ and $\beta$ are strongly Cartesian morphisms,
the functions $\alpha_{*}$ and $\beta_{*}$ are bijections. Therefore,
by the commutativity of the above diagram, the function
$(\alpha \circ \beta)_{*}$ is also a bijection and, hence,
$\alpha \circ \beta$ is strongly Cartesian. \qed

Taking $g = \id{T}$ in Definition \ref{def:stronglycartesian}, we see
that every strongly Cartesian morphism is Cartesian. The converse also
is true in fibred categories.

\begin{proposition}
  \label{pro:stronglycartfibred}
  (See \cite[Expos\'{e} VI, Proposition 6.11]{Grothendieck:1971:REG}.)
  Let $\mathcal{E}$ be a category, and $(\mathcal{F}, p)$ an
  $\mathcal{E}$-category. Then, $\mathcal{F}$ is a fibred category
  over $\mathcal{E}$ if and only if it satisfies the condition
  $\mathbf{Fib_I}$ of Definition \ref{def:fibredcategory}, and the
  following condition:
  \begin{enumerate}
  \item[$\mathbf{Fib_{II}'}$] Every Cartesian morphism in
    $\mathcal{F}$ is strongly Cartesian.
  \end{enumerate}
\end{proposition}

\noindent \emph{Proof}. Suppose that $\mathcal{F}$ satisfies condition
$\mathbf{Fib_I}$ of Definition \ref{def:fibredcategory}. Then, we need
to show that $\mathcal{F}$ satisfies $\mathbf{Fib_{II}}$ if and only
if it satisfies $\mathbf{Fib_{II}'}$.

So, let us suppose, first, that $\mathcal{F}$ satisfies
$\mathbf{Fib_{II}}$. We will prove that every Cartesian morphism
$\alpha : \eta \rightarrow \xi$ in $\mathcal{F}$ is strongly
Cartesian. Let $f : T \rightarrow S$, $g : U \rightarrow T$ and
$\zeta \in \Ob{\mathcal{F}_U}$ be as in Definition
\ref{def:stronglycartesian}. Let $\beta : \zeta' \rightarrow \eta$
be the inverse image of $\eta$ by $g$. Then, the diagram
\begin{displaymath}
  \xymatrix{
    \Hom[U]{\zeta}{\zeta'} \ar[rrr]^{\beta_{*}}
    \ar[ddrrr]_{(\alpha \circ \beta)_{*}} & & &
    \Hom[g]{\zeta}{\eta} \ar[dd]^{\alpha_{*}} \\
    & & & \\
    & & & \Hom[f \circ g]{\zeta}{\xi}
  }
\end{displaymath}
commutes. Since $(\zeta', \beta)$ is the inverse image of an
object, the morphism $\beta$ is Cartesian and, hence, the function
$\beta_{*}$ is a bijection. On the other hand, by
$\mathbf{Fib_{II}}$, the composition
$\alpha \circ \beta : \zeta' \rightarrow \xi$ is a Cartesian
morphism. Therefore, the function $(\alpha \circ \beta)_{*}$ is also
a bijection. Thus, by the commutativity of the above diagram, the
function $\alpha_{*}$ is a bijection, as well. In other words,
$\alpha$ is a strongly Cartesian morphism.

Conversely, if $\mathcal{F}$ satisfies $\mathbf{Fib_{II}'}$, then
every composition of Cartesian morphisms is a product of strongly
Cartesian morphisms, which, by Proposition
\ref{pro:stronglycartcompose}, is strongly Cartesian, hence
Cartesian.  Therefore, $\mathcal{F}$ satisfies the condition
$\mathbf{Fib_{II}}$. \qed

\section{Split categories}
\label{sec:SplitCategories}

We will now discuss a special class of fibred categories, which is
relevant to our description of binary graphs. We refer the reader to
\cite[Expos\'{e} VI, Sections 7, 8 and 9]{Grothendieck:1971:REG} for a
discussion, in full generality, of the topics mentioned in this
Section.

\begin{definition}
  \label{def:cleavage}
  (See
  \cite[Expos\'{e} VI, D\'{e}finition 7.1]{Grothendieck:1971:REG}.)
  Let $\mathcal{F}$ be an $\mathcal{E}$-category. A \emph{cleavage} of
  $\mathcal{F}$ over $\mathcal{E}$ is an attachment, to each morphism
  $f : T \rightarrow S$ in $\mathcal{E}$, of an inverse image functor
  by $f$ in $\mathcal{F}$, say, $f^{*}$. A \emph{cleaved category}
  over $\mathcal{E}$ is an $\mathcal{E}$-category $\mathcal{F}$,
  together with a cleavage of $\mathcal{F}$.
\end{definition}

\begin{definition}
  \label{def:splitting}
  (See \cite[Expos\'{e} VI, Section 9]{Grothendieck:1971:REG}, and
  \cite[Definition 4.1]{Streicher:1999:FCJ}.) Let $\mathcal{F}$ be an
  $\mathcal{E}$-category. We say that a cleavage of $\mathcal{F}$ over
  $\mathcal{E}$ is a \emph{splitting} if it satisfies the following
  conditions:
  \begin{enumerate}
  \item For every $S \in \Ob{\mathcal{E}}$, the inverse image functor
    $(\id{S})^{*} : \mathcal{F}_S \rightarrow \mathcal{F}_S$ equals
    the identity functor $\id{\mathcal{F}_S}$, and for every object
    $\xi \in \Ob{\mathcal{F}_S}$, we have
    $\alpha_{\id{S}} (\xi) = \id{\xi}$, where
    $\alpha_{\id{S}} (\xi) : (\id{S})^{*} \xi \rightarrow \xi$ is the
    canonical morphism.
  \item For every pair of composable morphisms, say,
    $f : T \rightarrow S$ and $g : U \rightarrow T$ in $\mathcal{E}$,
    we have
    \begin{gather*}
      (f \circ g)^{*} = g^{*} \circ f^{*}
      : \mathcal{F}_S \rightarrow \mathcal{F}_U, \quad \\
      \intertext{and for all $\xi \in \Ob{\mathcal{F}_S}$, we have}
      \alpha_{f \circ g} (\xi) = %
      \alpha_f (\xi) \circ \alpha_g (f^{*} \xi) : %
      (f \circ g)^{*} \xi = g^{*} f^{*} \xi \rightarrow \xi,
    \end{gather*}
    where $\alpha_f (\xi) : f^{*} \xi \rightarrow \xi$ is, as usual,
    the canonical morphism.
  \end{enumerate}
  A \emph{split category} over $\mathcal{E}$ is an
  $\mathcal{E}$-category $\mathcal{F}$, together with a splitting of
  $\mathcal{F}$. If $\mathcal{F}$ and $\mathcal{G}$ are split
  categories over $\mathcal{E}$, then, a
  \emph{morphism of split categories}, from $\mathcal{F}$ to
  $\mathcal{G}$, is an $\mathcal{E}$-functor
  $F : \mathcal{F} \rightarrow \mathcal{G}$, such that for every
  morphism $f : T \rightarrow S$ in $\mathcal{E}$, and for every
  $\xi \in \Ob{\mathcal{F}_S}$, we have
  $F (f^{*} \xi) = f^{*} F (\xi)$ and
  $F (\alpha_f (\xi)) = \alpha_f (F (\xi))$. We thus obtain a category
  $\Split{\mathcal{E}}$ of split categories over $\mathcal{E}$.
\end{definition}

\begin{notation}
  \label{not:CatHomCat}
  Let $\Cat$ denote the ``category of categories'', i.e., the category
  whose objects are categories, and whose morphisms are functors
  between categories. For any pair of categories $\mathcal{A}$ and
  $\mathcal{B}$, let $\HomCat{\mathcal{A}}{\mathcal{B}}$ denote the
  category whose objects are functors from $\mathcal{A}$ to
  $\mathcal{B}$, and whose morphisms are natural transformations of
  functors. In particular,
  $\HomCat{\mathcal{A}^\circ}{\mathcal{B}}$,
  where $\mathcal{A}^\circ$ is the opposite category of $\mathcal{A}$,
  is the category of contravariant functors from $\mathcal{A}$ to
  $\mathcal{B}$.
\end{notation}

\begin{remark}
  \label{rem:SplitHomEquivalence}
  Given a split category $\mathcal{F}$ over $\mathcal{E}$, we obtain a
  functor
  $\phi (\mathcal{F}) : \mathcal{E}^\circ \rightarrow \Cat$, that is
  defined as follows:
  \begin{gather*}
    \phi (\mathcal{F}) : \Ob{\mathcal{E}^\circ} \rightarrow \Ob{\Cat},
    \quad S \mapsto \mathcal{F}_S, \quad \text{and} \\
    \phi (\mathcal{F}) : \Hom{T}{S} \rightarrow
    \Hom{\mathcal{F}_S}{\mathcal{F}_T}, \quad f \mapsto f^{*}.
  \end{gather*}
  If $F : \mathcal{F} \rightarrow \mathcal{G}$ is a morphism in
  $\Split{\mathcal{E}}$, we obtain a natural transformation
  $\phi (F) : \phi (\mathcal{F}) \rightarrow \phi (\mathcal{G})$ as
  follows: given any object $S \in \Ob{\mathcal{E}}$, we define
  \begin{displaymath}
    \phi (F)_S : \phi (\mathcal{F}) (S) = \mathcal{F}_S \rightarrow
    \phi (\mathcal{G}) (S) = \mathcal{G}_S, \quad
    \phi (F)_S = \restrict{F}{\mathcal{F}_S}.
  \end{displaymath}
  We thus obtain a functor
  \begin{displaymath}
    \phi : \Split{\mathcal{E}} \rightarrow
    \HomCat{\mathcal{E}^\circ}{\Cat}.  
  \end{displaymath}
  This functor, $\phi$, is an equivalence of categories (see
  \cite[Expos\'{e} VI, Section 9]{Grothendieck:1971:REG}).
\end{remark}

\begin{proposition}
  \label{pro:SplitFibred}
  Every split category $\mathcal{F}$, over a category $\mathcal{E}$,
  is a fibred category.
\end{proposition}

\noindent \emph{Proof}. Since $\mathcal{F}$ admits a cleavage, it
satisfies the condition $\mathbf{Fib_I}$ of Definition
\ref{def:fibredcategory}. We will show that $\mathcal{F}$ satisfies
the condition $\mathbf{Fib_{II}'}$ of Proposition
\ref{pro:stronglycartfibred}, i.e., that every Cartesian morphism in
$\mathcal{F}$ is strongly Cartesian. Let
$\alpha : \eta \rightarrow \xi$ be a Cartesian morphism in
$\mathcal{F}$, and let $f : T \rightarrow S$ denote the morphism
$p (f)$ in $\mathcal{E}$, where
$p : \mathcal{F} \rightarrow \mathcal{E}$ is the canonical functor.
Let $g : U \rightarrow T$ be a morphism in $\mathcal{E}$, and let
$\zeta \in \Ob{\mathcal{F}_U}$. We, then, have to show that the
function
\begin{align*}
  \alpha_{*} : \Hom[g]{\zeta}{\eta} & \rightarrow
  \Hom[f \circ g]{\zeta}{\xi}, \\
  v & \mapsto \alpha \circ v,
\end{align*}
is a bijection. Since $\alpha$ is Cartesian, the pair $(\eta, \alpha)$
is an inverse image of $\xi$ by $f$. On the other hand, the given
splitting of $\mathcal{F}$ provides another inverse image
$(f^{*} \xi, \alpha_f (\xi))$ of $\xi$ by $f$. Therefore, by the
definition of Cartesian morphisms, there exists a unique
$T$-isomorphism $\beta : f^{*} \xi \rightarrow \eta$ such that
$\alpha_f (\xi) = \alpha \circ \beta$. We now have the following
commutative diagram:
\begin{displaymath}
  \xymatrix{%
    \Hom[g]{\zeta}{\eta} \ar[rr]^{\alpha_{*}} & &
    \Hom[f \circ g]{\zeta}{\xi} \\
    & & \\
    \Hom[g]{\zeta}{f^{*} \xi} \ar[uu]^{\beta_{*}} & &
    \Hom[U]{\zeta}{(f \circ g)^{*} \xi}
    \ar[ll]^{\alpha_g (f^{*} \xi)_{*}}
    \ar[uu]_{\alpha_{f \circ g} (\xi)_{*}} \\
    & & \\
    & \Hom[U]{\zeta}{g^{*} f^{*} \xi}
    \ar[luu]_{\alpha_g (f^{*} \xi)_{*}}
    \ar@{=}[ruu] &
  }
\end{displaymath}
The triangle in the bottom of the diagram is defined because
$g^{*} f^{*} \xi = (f \circ g)^{*} \xi$, and the commutativity of the
rectangle follows from the equations
\begin{displaymath}
  \alpha \circ \beta \circ \alpha_g (f^{*} \xi) =
  \alpha_f (\xi) \circ \alpha_g (f^{*} \xi) =
  \alpha_{f \circ g} (\xi).
\end{displaymath}
The function $\beta_{*}$ is a bijection because $\beta$ is an
isomorphism, while the bottom and right sides of the rectangle are
bijections because the morphisms $\alpha_g (f^{*} \xi)$ and
$\alpha_{f \circ g} (\xi)$, respectively, are Cartesian. Therefore,
the the top edge of the rectangle, too, is a bijection. This proves
that $\alpha : \eta \rightarrow \xi$ is strongly Cartesian. \qed

\begin{definition}
  \label{def:RigidReduced}
  (See \cite[Expos\'{e} VI, Section 9]{Grothendieck:1971:REG}.) We say
  that a category $\mathcal{A}$ is \emph{rigid} if, for every object
  $\xi \in \Ob{\mathcal{A}}$, the identity morphism $\id{\xi}$ is the
  only automorphism of $\xi$. The category $\mathcal{A}$ is said to be
  \emph{reduced}, if whenever two objects in $\mathcal{A}$ are
  isomorphic, they are, in fact, equal.
\end{definition}

\begin{proposition}
  \label{pro:UniqueCleavage}
  Let $\mathcal{E}$ be a category, and let $\mathcal{F}$ be an
  $\mathcal{E}$-category. Suppose that for every object
  $S \in \Ob{\mathcal{E}}$, the categorical fibre $\mathcal{F}_S$ is a
  rigid and reduced category. Then, $\mathcal{F}$ is a fibred category
  if and only if it is split. If this is, indeed, the case, then there
  exists a unique cleavage of $\mathcal{F}$ over $\mathcal{E}$, and
  that cleavage is a splitting.
\end{proposition}

\noindent \emph{Proof}. By Proposition \ref{pro:SplitFibred}, if
$\mathcal{F}$ is split over $\mathcal{E}$, then it is a fibred
category. Conversely, suppose that $\mathcal{F}$ is a fibred category
over $\mathcal{E}$. We will, then, show that $\mathcal{F}$ has a
unique cleavage, and that this cleavage is a splitting of
$\mathcal{F}$ over $\mathcal{E}$.

Since $\mathcal{F}$ is a fibred category over $\mathcal{E}$, using the
Axiom of Choice, as in \cite[Section 4]{Streicher:1999:FCJ}, we see
that $\mathcal{F}$ admits a cleavage. Suppose that $\mathcal{F}$
admits two cleavages $*$ and $\dagger$. For any morphism
$f : T \rightarrow S$ in $\mathcal{E}$, and for any object
$\xi \in \Ob{\mathcal{F}_S}$, denote the inverse images of $\xi$ by
$f$, with respect to the cleavages $*$ and $\dagger$, by
$(f^{*} \xi, \alpha_f (\xi))$ and $(f^{\dagger} \xi, \beta_f (\xi))$,
respectively. Then, there exists a unique $T$-isomorphism
$\theta_f (\xi) : f^{*} \xi \rightarrow f^{\dagger} \xi$, such that
$\alpha_f (\xi) = \beta_f (\xi) \circ \theta_f (\xi)$. Since
$\mathcal{F}_T$ is a reduced category, the $T$-isomorphic objects
$f^{*} \xi$ and $f^{\dagger} \xi$ are equal. Now, because
$\mathcal{F}_T$ is a rigid category, the $T$-automorphism
$\theta_f (\xi)$ of $f^{*} \xi$ equals the identity morphism.
Therefore, $\alpha_f (\xi) = \beta_f (\xi)$. In other words, for all
$f$ and $\xi$, the inverse images $(f^{*} \xi, \alpha_f (\xi))$ and
$(f^{\dagger} \xi, \beta_f (\xi))$ are equal. Thus, the two cleavages
$*$ and $\dagger$ are equal and, so, $\mathcal{F}$ has a unique
cleavage over $\mathcal{E}$.

We will show, next, that the unique cleavage of $\mathcal{F}$ is a
splitting. Let $f$ and $\xi$ be as above, and let
$g : U \rightarrow T$ be another morphism in $\mathcal{E}$. Then, by
the condition $\mathbf{Fib_{II}}$ of Definition
\ref{def:fibredcategory}, the composition
$\alpha_f (\xi) \circ \alpha_g (f^{*} \xi) : g^{*} f^{*} \xi
\rightarrow \xi$
is a Cartesian morphism over $f \circ g : U \rightarrow S$.
Therefore, by the universal property of inverse images, there exists a
unique $U$-isomorphism
$c_{f, g} (\xi) : g^{*} f^{*} \xi \rightarrow (f \circ g)^{*} \xi$,
such that the following diagram commutes:
\begin{displaymath}
  \xymatrix{%
    g^{*} f^{*} \xi \ar[rr]^{\alpha_g (f^{*} \xi)}
    \ar@{-->}[dd]_{c_{f, g} (\xi)} & &
    f^{*} \xi \ar[dd]^{\alpha_f (\xi)} \\
    & & \\
    (f \circ g)^{*} \xi \ar[rr]_{\alpha_{f \circ g} (\xi)}
    & & \xi
  }
\end{displaymath}
Because $\mathcal{F}_U$ is a rigid and reduced category, we know that
the $U$-isomorphic objects $g^{*} f^{*} \xi$ and $(f \circ g)^{*} \xi$
are equal, and that the $U$-automorphism $c_{f, g} (\xi)$ of
$(f \circ g)^{*} \xi$ is the identity morphism. Therefore, we get
$\alpha_{f \circ g} (\xi) = \alpha_f (\xi) \circ \alpha_g (f^{*}
\xi)$.
In a similar manner, since $(\xi, \id{\xi})$ and
$((\id{S})^{*} \xi, \alpha_{\id{S}} (\xi))$ are both inverse images of
$\xi$ by $\id{S}$, we see that
$((\id{S})^{*} \xi, \alpha_{\id{S}} (\xi)) = (\xi, \id{\xi})$. We
conclude that the unique cleavage of $\mathcal{F}$ is a splitting over
$\mathcal{E}$. \qed

\section{Fibred category of binary automata}
\label{sec:binaryfibredcat}

Recall that in Section \ref{sec:binaryautomata}, we have defined a
functor $p : \mathcal{BA} \rightarrow \mathcal{BG}_{*}$, from the
category $\mathcal{BA}$ of binary automata to the category
$\mathcal{BG}_{*}$ of pointed binary graphs. In the terminology of
Section \ref{sec:fibredcategories}, $(\mathcal{BA}, p)$ is a
$\mathcal{BG}_{*}$-category. We will now show that it is, in fact, a
fibred category over $\mathcal{BG}_{*}$. But, first, let us note that
we have a canonical functor in the other direction, too.

Given a pointed binary graph $P = (\Gamma, q_0)$, where
$\Gamma = (Q, \delta)$ we get an admissible binary automaton
$M_{P} = (Q, \Sigma, \delta, q_0, F)$, where $F$ is the set of leaves
in $Q$. If $P = (\Gamma, q_0)$, and $P' = (\Gamma', q_0')$ are two
pointed binary graphs, and if $f : P \rightarrow P'$ is a morphism,
then the underlying function $f : Q \rightarrow Q'$, being
$\Sigma^{*}$-equivariant, takes the leaves of $\Gamma$ to leaves of
$\Gamma'$ and, so, defines a morphism of binary automata $M_f : M_P
\rightarrow M_{P'}$. We thus obtain a functor $s$ from the category
$\mathcal{BG}_{*}$ to the category $\mathcal{BA}$,
\begin{gather*}
  s: \Ob{\mathcal{BG}_{*}} \rightarrow \Ob{\mathcal{BA}}, \quad 
  s (P) = M_P, \quad \text{and} \\
  s : \Hom{P}{P'} \rightarrow \Hom{M_P}{M_{P'}}, \quad s (f) = M_f.
\end{gather*}

\begin{proposition}
  \label{pro:sectionfibred}
  The above functor $s: \mathcal{BG}_{*} \rightarrow \mathcal{BA}$ is
  a section of the $\mathcal{BG}_{*}$-category $(\mathcal{BA}, p)$.
  This section is an isomorphism of categories from
  $\mathcal{BG}_{*}$, to the full subcategory of $\mathcal{BA}$
  consisting of admissible automata.
\end{proposition}

\noindent \emph{Proof}. It is clear that
$p \circ s = \id{\mathcal{BG}_{*}}$. From the definition of
admissibility, it follows that for every pointed binary graph $P$, the
binary automaton $M_P$ is admissible. Let
$M = (Q, \Sigma, \delta, q_0, F)$ be an admissible automaton, and let
$P_M = (\Gamma_M, q_0)$. Then, it follows, from Remark
\ref{rem:admissible}, that the set of leaves of $\Gamma_M$ equals the
set, $F$, of final states of $M$. Thus $s (P_M) = M$, i.e.,
$s \circ p = \id{}$ on $\Ob{\mathcal{ABA}}$, where $\mathcal{ABA}$ is
the full subcategory of $\mathcal{BA}$, consisting of admissible
binary automata. Moreover, if $f : M \rightarrow M'$ is a morphism in
$\mathcal{ABA}$, then
$s \circ p (f) = s (P_f) = M_{P_f} = f$. Thus, $s \circ p = \id{}$ on
$\Hom{M}{M'}$, and hence $s \circ p = \id{\mathcal{ABA}}$. \qed

We, thus, obtain an identification between pointed binary graphs and
admissible binary automata.

\begin{remark}
  \label{rem:binaryfibre}
  It follows, from Proposition \ref{pro:sectionfibred}, that
  $s : \Ob{\mathcal{BG}_{*}} \rightarrow \Ob{\mathcal{BA}}$ is an
  injective function. On the other hand,
  $p : \Ob{\mathcal{BA}} \rightarrow \Ob{\mathcal{BG}_{*}}$ is a
  many-to-one function. Indeed if
  $P = (\Gamma, q_0) \in \Ob{\mathcal{BG}_{*}}$, and if
  $\Gamma = (Q, \delta)$, then for every subset $F \subset Q$, we
  obtain an object $M = (Q, \Sigma, \delta, q_0, F)$ such that
  $p (M) = P$. Indeed, the categorical fibre, $\mathcal{BA}_P$, of
  $\mathcal{BA}$ over $P$ can be described as follows:
  \begin{enumerate}
  \item $\Ob{\mathcal{BA}_P} = \powerset{Q}$, the power set of $Q$.
  \item If $F_1$ and $F_2$ are subsets of $Q$, then
    $\Hom[P]{F_1}{F_2}$ is the empty set if $F_1$ is not a subset of
    $F_2$, and is the singleton $\set{\imath_{F_1, F_2}}$,
    if $F_1 \subset F_2$, in which case $\imath_{F_1, F_2}$ is the
    inclusion map $F_1 \hookrightarrow F_2$.
  \end{enumerate}
  We call this category the
  \emph{posetal category defined by the power set of Q}, and denote it
  by $\mathcal{P}_Q$.
\end{remark}

\begin{definition}
  \label{def:BinaryCleavage}
  Let $P = (\Gamma, q_0)$, and $P' = (\Gamma', q_0')$ be two pointed
  binary graphs, where $\Gamma = (Q, \delta)$ and
  $\Gamma' = (Q', \delta')$, and let $f : P \rightarrow P'$ be a
  morphism. Let $M' \in \Ob{\mathcal{BA}_{P'}}$ be a binary automaton
  over $P'$, say, $M' = (Q', \Sigma, \delta', q_0', F')$. We, then,
  define a binary automaton $f^{*} M' \in \Ob{\mathcal{BA}_P}$ by
  $f^{*} M' = (Q, \Sigma, \delta, q_0, f^{-1} (F'))$. Further, we
  define a morphism of binary automata,
  $\alpha_f (M') : f^{*} M' \rightarrow M'$ by setting
  $\alpha_f (M') = f : Q \rightarrow Q'$. We will show below (see the
  proof of Theorem \ref{thm:BinarySplitCategory}), that the pair
  $(f^{*} M', \alpha_f (M'))$ is an inverse image of $M'$ by $f$. We
  thus obtain a cleavage of $\mathcal{BA}$ over $\mathcal{BG}_{*}$. We
  call it the \emph{canonical cleavage} of $\mathcal{BA}$ over
  $\mathcal{BG}_{*}$.
\end{definition}

\begin{theorem}
  \label{thm:BinarySplitCategory}
  The $\mathcal{BG}_{*}$-category $(\mathcal{BA}, p)$ is a fibred
  category, in fact, a split category. The canonical cleavage of
  $\mathcal{BA}$ (see Definition \ref{def:BinaryCleavage}) is its only
  cleavage, and this cleavage is a splitting.
\end{theorem}

\noindent \emph{Proof}. We will, first, prove that the canonical
cleavage is, indeed, a cleavage of $\mathcal{BA}$ over
$\mathcal{BG}_{*}$. Let $P = (\Gamma, q_0)$, and
$P' = (\Gamma', q_0')$ be two pointed binary graphs, where $\Gamma =
(Q, \delta)$ and $\Gamma' = (Q', \delta')$, and let $f : P \rightarrow
P'$ be a morphism. Let $M' \in \Ob{\mathcal{BA}_{P'}}$ be a binary
automaton over $P'$, say, $M' = (Q', \Sigma, \delta', q_0', F')$. Let
$(f^{*} M', \alpha_f (M'))$ be as in Definition
\ref{def:BinaryCleavage}. Let $M'' \in \Ob{\mathcal{BA}_P}$ be a
binary automaton over $P$, say, $M'' = (Q, \Sigma, \delta, q_0, F'')$,
and let $u : M'' \rightarrow M'$ be an $f$-morphism of binary
automata. We need to show that there exists a unique $P$-morphism
$\overline{u} : M'' \rightarrow f^{*} M'$ such that the diagram
\begin{displaymath}
  \xymatrix{%
    M'' \ar[drrr]^u \ar@{-->}[dr]_{\overline{u}} & & & \\
    & f^{*} M' \ar[rr]_{\alpha} & & M'
  }
\end{displaymath}
commutes. Since $u$ is an $f$-morphism, by the definition of the
functor $p : \mathcal{BA} \rightarrow \mathcal{BG}_{*}$, we see that
the underlying function $u : Q \rightarrow Q'$ of $u$, equals $f$.
Further, since $u$ is a morphism of binary automata, we have
$f (F'') = u (F'') \subset F'$, i.e., we have
$F'' \subset f^{-1} (F')$. Now, define
$\overline{u} : M'' \rightarrow f^{*} M'$, by
$\overline{u} = \id{Q} : Q \rightarrow Q$. Then, the fact that $F''$
is a subset of $f^{-1} (F')$ implies that $\overline{u}$ is a morphism
of binary automata. It is, obviously, a $P$-morphism and, moreover,
\begin{displaymath}
  \alpha_f (M') \circ \overline{u} = f \circ \id{Q} = f = u.
\end{displaymath}
This proves the ``existence'' part of Cartesianness. To prove the
``uniqueness'' part, suppose that $v$ and $w$ are $P$-morphisms from
$M''$ to $f^{*} M'$, such that
$\alpha_f (M') \circ v = \alpha_f (M') \circ w$. Since
$p (v) = \id{P}$, the underlying function $v : Q \rightarrow Q$ equals
$\id{Q}$. Similarly, $w$ equals $\id{Q}$, hence $v = w$. We have,
thus, proved that the canonical cleavage is, in fact, a cleavage of
$\mathcal{BA}$ over $\mathcal{BG}_{*}$.

We will, next, show that the canonical cleavage of $\mathcal{BA}$ is a
splitting. If $f = \id{P'}$, then $f^{-1} (F') = F'$, hence
$(\id{P'})^{*} M' = M'$, and $\alpha_{\id{P'}} (M') = \id{M'}$. Next,
suppose that $f$ and $M'$ are as above, and let $g : P'' \rightarrow P$
be another morphism in $\mathcal{BG}_{*}$. Then,
$(f \circ g)^{-1} (F') = g^{-1} (f^{-1} (M'))$, hence
$(f \circ g)^{*} M' = g^{*} f^{*} M'$. Further, since the underlying
function of a morphism of the type $\alpha_f (M')$, as we saw above,
equals $f$, we have
\begin{displaymath}
  \alpha_{f \circ g} (M') = f \circ g =
  \alpha_f (M') \circ \alpha_g (f^{*} M').
\end{displaymath}
Therefore, the canonical cleavage of $\mathcal{BA}$ is a splitting
over $\mathcal{BG}_{*}$. Thus, $(\mathcal{BA}, p)$ is a split category
and, so, by Proposition \ref{pro:SplitFibred}, it is a fibred category
over $\mathcal{BG}_{*}$.

Finally, we have observed (see Remark \ref{rem:binaryfibre}), that for
every object $P \in \Ob{\mathcal{BG}_{*}}$, the categorical fibre
$\mathcal{BA}_P$ is the posetal category $\powersetcat{Q}$ defined by
the power set of $Q$, where $Q$ is the set of nodes of $P$. Thus,
every categorical fibre of $\mathcal{BA}$ is a rigid and reduced
category. Therefore, by Proposition \ref{pro:UniqueCleavage}, the
fibred category $\mathcal{BA}$ has a unique cleavage. \qed

\begin{remark}
  \label{rem:BinaryFunctor}
  It follows from Remark \ref{rem:SplitHomEquivalence} that there is a
  canonical equivalence $\phi$ from the category
  $\Split{\mathcal{BG}_{*}}$ to the category
  $\HomCat{\mathcal{BG}_{*}^{\circ}}{\Cat}$. From the definition of
  $\phi$, we see that the functor
  $\phi (\mathcal{BA}) : \mathcal{BG}_{*}^{\circ} \rightarrow \Cat$
  is isomorphic to the following functor:
  \begin{align*}
    \Ob{\mathcal{BG}_{*}^{\circ}} & \rightarrow \Cat, \quad %
    P \mapsto \powersetcat{Q}, \quad \text{and} \\
    \Hom{P}{P'} & \rightarrow \Hom{\powersetcat{Q'}}{\powersetcat{Q}},
    \quad f \mapsto f^{*}.
  \end{align*}
  We use here the notation $Q$ (respectively, $Q'$) for the set of
  nodes of the pointed binary graph $P$ (respectively, $P'$); if $Q$
  is a set, $\powersetcat{Q}$ denotes the posetal category defined by
  the power set of $Q$ (see Remark \ref{rem:binaryfibre}); and, for
  any function $f : Q \rightarrow Q'$, the symbol $f^{*}$ denotes the
  obvious pull-back functor from $\powersetcat{Q'}$ to
  $\powersetcat{Q}$.
\end{remark}

\section{Conclusion}
\label{sec:conclusion}

In this article, we have developed a set-theoretic formalism for
binary graphs. We have expressed various notions regarding graphs in
terms of the inheritance order on the set of nodes, and in terms of
the action of bit strings on nodes. We have exhibited pointed binary
graphs as a subcategory of transition systems. Together with the
set-theoretic notion of binary graphs, another interesting result, for
us, in this paper is the fact that binary automata form a fibred
category over pointed binary graphs. We would like to have a better
understanding of this fibred category.

We suggest that it would be interesting to formulate various
algorithms for trees, using the formalism developed here.  Obviously,
this formalism, alone, would not affect the efficiency of such
algorithms. However, we feel that the abstract set-theoretic
descriptions provided here would help in specifying the algorithms in
a clear and precise manner.

\section*{Acknowledgements}
\label{sec:acknowledgements}

I thank M.~Vidyasagar for providing me an opportunity to interact with
computational biologists, and Sharmila Mande and Ganesh Gayatri for
introducing me to the algorithms that motivated this work. I am
grateful to K.~Viswanath for his comments on an earlier version of
this article. I thank I.B.S.~Passi for asking me to lecture in his
seminar on geometric group theory, at the Harish-Chandra Research
Institute, Allahabad; the talks I gave there, on automatic groups,
motivated me to study finite automata.

\bibliographystyle{raghu-amsalpha}
\nocite{Abramsky:1995:HLC}
\nocite{Winskel:1995:MC}
\bibliography{raghu-strings,raghu-refs,raghu-crossrefs}

\ifx \path \undefined \input path.sty \fi
\providecommand{\eprintaux}[1]{\csname #1URL\endcsname}
\providecommand{\eprint}[3]{\href{\eprintaux{#1}#2}{\path|#3|}}
\providecommand{\arXivURL}{http://www.arXiv.org/abs/}
\providecommand{\ResearchIndexURL}{http://citeseer.nj.nec.com/}
\providecommand{\DigitalMathArchiveURL}{%
  http://www.sunsite.ubc.ca/DigitalMathArchive/}
\providecommand{\biburl}[1]{URL \href{#1}{\path|#1|}}
\newenvironment{bibannote}%
  {\begin{quotation} \noindent \sffamily \small}%
  {\end{quotation}}
\newenvironment{bibrefer}%
  {\begin{quotation} \noindent \sffamily \small}%
  {\end{quotation}}
\providecommand{\allcaps}[1]{{#1}}\providecommand{\singleletter}[1]{#1}
\providecommand{\bysame}{\leavevmode\hbox to3em{\hrulefill}\thinspace}
\providecommand{\MR}{\relax\ifhmode\unskip\space\fi MR }
\providecommand{\MRhref}[2]{%
  \href{http://www.ams.org/mathscinet-getitem?mr=#1}{#2}
}
\providecommand{\href}[2]{#2}
\begin{thebibliography}{AGM95}

\bibitem[AGM95]{Abramsky:1995:HLC}
S.~Abramsky, D.M. Gabbay, and T.S.E Maibaum (eds.), \emph{Handbook of logic in
  computer science, \allcaps{V}ol. 4, \allcaps{S}emantic modelling}, Oxford
  University Press, Oxford, UK, 1995, ISBN 0-19-853780-8.



\bibitem[AHU74]{Aho:1974:DAC}
Alfred~V. Aho, John~E. Hopcroft, and Jeffrey~D. Ullman, \emph{The design and
  analysis of computer algorithms}, Addison-Wesley, Boston 1974, ISBN
  981-4053-19-8.



\bibitem[GR71]{Grothendieck:1971:REG}
Alexandre Grothendieck and Michele Raynaud, \emph{Rev{\^e}tements {\'e}tales et
  groupe fondamental (\allcaps{SGA}{~1}), \allcaps{U}n \allcaps{S}{\'e}minaire
  dirig{\'e} par \allcaps{A}. \allcaps{G}rothendieck, \allcaps{A}ugment{\'e} de
  deux expos{\'e}s par \allcaps{M}. \allcaps{R}aynaud, \allcaps{S}{\'e}minaire
  de \allcaps{G}{\'e}om{\'e}trie \allcaps{A}lg{\'e}brique du
  \allcaps{B}ois-\allcaps{M}arie, 1960-61}, Lecture Notes in Mathematics, no.
  224, Springer-Verlag, Berlin, Heidelberg and New York, 1971,
  \eprint{arXiv}{math.AG/0206203}{arXiv:math.AG/0206203}, ISBN 3-540-05614-9.



\bibitem[HU79]{Hopcroft:1979:IAT}
John~E. Hopcroft and Jeffrey~D. Ullman, \emph{Introduction to automata theory,
  languages, and computation}, first ed., Addison-Wesley, Boston, 1979, ISBN
  0-201-02988-X.



\bibitem[Knu97]{Knuth:1997:FA}
Donald~E. Knuth, \emph{The art of computer programming, \allcaps{V}ol. 1,
  \allcaps{F}undamental algorithms}, third ed., Addison-Wesley, Boston, 1997,
  ISBN 81-7808-111-3.



\bibitem[Str99]{Streicher:1999:FCJ}
Thomas Streicher, \emph{Fibred categories {\`a} la \allcaps{J}.
  \allcaps{B}{\'e}nabou}, Notes of a course given at a spring school in Munich,
  \biburl{http://www.mathematik.tu-darmstadt.de/~streicher/}, 1999.



\bibitem[WN93]{Winskel:1993:MC}
Glynn Winskel and Mogens Nielsen, \emph{Models for concurrency}, Tech. Report
  PB-463, Department of Computer Science, University of Aarhus, Aarhus,
  Denmark, November 1993, \biburl{http://www.daimi.au.dk/PB/463/}, Preprint
  version of {\cite{Winskel:1995:MC}}.



\bibitem[WN95]{Winskel:1995:MC}
\bysame, \emph{Models for concurrency}, in Abramsky et~al.
  \cite{Abramsky:1995:HLC}, pp.~1--148.



\end{thebibliography}

\end{document}